\documentclass[12pt]{amsart}
\usepackage{graphicx,amsmath,amssymb}

\numberwithin{equation}{section}

\newcommand{\C}{\mathbb C}
\newcommand{\R}{\mathbb R}
\newcommand{\Z}{\mathbb Z}
\newcommand{\N}{\mathbb N}
\renewcommand{\d}{\prime} 
\newcommand{\dd}{{\prime \prime}}

\renewcommand{\Re}{{\rm Re}\,}

\newtheorem{theorem}{Theorem}[section]
\newtheorem{lemma}[theorem]{Lemma}
 
\newtheorem{corollary}[theorem]{Corollary}

\newtheorem{definition}{Definition}
\newtheorem*{remark}{Remark}
\newtheorem*{remarks}{Remarks}

\begin{document}
\title[]
{Schr\"odinger type eigenvalue problems  with polynomial potentials: Asymptotics of  eigenvalues}
\author[]
{Kwang C. Shin}
\address{Department of Mathematics, University of Missouri, Columbia, MO 65211, USA}
\date{November 4, 2004}

\begin{abstract}
For integers $m\geq 3$ and $1\leq\ell\leq m-1$, we study the eigenvalue problem $-u^\dd(z)+[(-1)^{\ell}(iz)^m-P(iz)]u(z)=\lambda u(z)$ with the boundary conditions that $u(z)$ decays to zero as $z$ tends to infinity along the rays $\arg z=-\frac{\pi}{2}\pm \frac{(\ell+1)\pi}{m+2}$ in the complex plane, where $P(z)=a_1 z^{m-1}+a_2 z^{m-2}+\cdots+a_{m-1} z$ is a polynomial. We provide asymptotic expansions of the eigenvalue counting function and the eigenvalues $\lambda_{n}$. Then we apply these to the inverse spectral problem, reconstructing some coefficients of polynomial potentials from asymptotic expansions of the eigenvalues. Also, we show for arbitrary $\mathcal{PT}$-symmetric polynomial potentials of degree $m\geq 3$ and all symmetric decaying boundary conditions that the eigenvalues are all real and positive, with only finitely many exceptions.
\end{abstract}

\maketitle

\begin{center}
{\it Preprint.}
%{\it To appear in: When the paper has been accepted, the Journal Title goes here}
\end{center}

\baselineskip = 18pt

\section{Introduction}
\label{introduction}
In this paper, we study  Schr\"odinger operators with any polynomial potential of degree $m\geq 3$ with complex coefficients, under   decaying boundary conditions along two rays to infinity in the complex plane, and provide asymptotic expansions of the eigenvalue counting functions and the eigenvalues. Then we will use these to  reconstruct some coefficients of polynomial potentials from asymptotic expansions of the eigenvalues,  and to show that all but finitely many eigenvalues of every $\mathcal{PT}$-symmetric oscillator with  a polynomial potential are real and positive.

For integers $m\geq 3$ fixed and $1\leq \ell\leq m-1$, we are considering the non-standard eigenvalue problems
\begin{equation}\label{ptsym}
H_{\ell,P} u(z,\lambda):=\left[-\frac{d^2}{dz^2}+(-1)^{\ell}(iz)^m-P(iz)\right]u(z,\lambda)=\lambda u(z,\lambda),\quad\text{for some $\lambda\in\C$},
\end{equation}
with the boundary condition that 
\begin{equation}\label{bdcond}
\text{$u(z,\lambda)\rightarrow 0$ exponentially, as $z\rightarrow \infty$ along the two rays}\quad \arg z=-\frac{\pi}{2}\pm \frac{(\ell+1)\pi}{m+2},
\end{equation}
where $P$ is a polynomial of degree at most $m-1$ of the form 
\begin{equation}\nonumber
P(z)=a_1z^{m-1}+a_2z^{m-2}+\cdots+a_{m-1}z,\quad a_j\in\C\,\,\text{\,for $1\leq j\leq m-1$}.
\end{equation}

The anharmonic oscillators $H_{\ell,P}$ with the various boundary conditions \eqref{bdcond} are considered in \cite{Bender2,Shin2}. The most studied case is 
when $m$ is even and $\ell=\frac{m}{2}$, for example, see \cite{AA,Bender-1,Bender3,HR,HMM,MAS,Simon}. In this case, $H_{\ell,P}$ is an Schr\"odinger operator in $L^2(\R)$. This is self-adjoint if all the coefficients of the polynomial $(-1)^{\ell}(iz)^m-P(iz)$ in $z$ are real, and non-self-adjoint if a coefficient of $(-1)^{\ell}(iz)^m-P(iz)$ is non-real. The case when $\ell=1$ has been studied extensively in recent years in the context of $\mathcal{PT}$-symmetric theory \cite{Bender,Dorey,Handy2,Shin,Sibuya}.

Throughout this paper, we use the integer $m\geq 3$ for the degree of the polynomial potential and integer $\ell$ with  $1\leq\ell\leq m-1$. We will use $\lambda_n$, depending on the potential and the boundary condition, to denote  the eigenvalues of $H_{\ell,P}$, without explicitly indicating their dependence on the potential and the boundary condition. 
Also, we let $$a:=(a_1,a_2,\ldots, a_{m-1})\in\C^{m-1}$$ be the coefficient vector of $P(z)$. 

If a nonconstant function $u$ satisfies \eqref{ptsym} with some $\lambda\in\C$ and the boundary condition \eqref{bdcond}, then we call $\lambda$ an {\it eigenvalue} of $H_{\ell,P}$ and $u$ an {\it eigenfunction of $H_{\ell,P}$ associated with the eigenvalue $\lambda$}. Also, the {\it geometric multiplicity of an eigenvalue $\lambda$} is the number of linearly independent eigenfunctions  associated with the eigenvalue $\lambda$. For each $\ell$ fixed, we number the eigenvalues $\{\lambda_{n}\}_{n\geq n_0}$ in the order of nondecreasing magnitudes, counting their ``algebraic multiplicities'', where the integer $n_0$, depending on the potential and the boundary condition, is due to our method of proof of Theorem \ref{main_thm1}. In Theorem \ref{main_thm1} we show that for every large $n\in\N$, there exists $\lambda_n$ satisfying  \eqref{main_result} below. However, we do not know the number of eigenvalues ``near'' zero, and this is why we need to have $n_0$ in numbering the eigenvalues.

Before we state our main theorems, we first introduce some known facts by Sibuya \cite{Sibuya} about the eigenvalues $\lambda_n$ of $H_{\ell,P}$.
\begin{theorem}\label{main2}
The eigenvalues $\lambda_{n}$ of $H_{\ell,P}$ have the following properties.
\begin{enumerate}
\item[(I)] The set of all eigenvalues is a  discrete set in $\C$.
\item[(II)] The geometric multiplicity of every eigenvalue is one.
\item[(III)] Infinitely many eigenvalues, accumulating at infinity, exist.
\item[(IV)] When $\ell=1$ the eigenvalues have the following asymptotic expansion
\begin{equation}\label{Bender_exp}
\lambda_{n} =\left(\frac{\Gamma\left(\frac{3}{2}+\frac{1}{m}\right)\sqrt{\pi}\left(n+\frac{1}{2}\right)}{\sin \left(\frac{\pi}{m}\right)\Gamma\left(1+\frac{1}{m}\right) } \right)^{\frac{2m}{m+2}}[1+o(1)]\quad\text{as $n$ tends to infinity},\quad n \in \N.
\end{equation}
\end{enumerate} 
\end{theorem}

This paper contains results on direct and inverse spectral probelms, and their applications to $\mathcal{PT}$-symmetric potential problems. Theorem \ref{main_thm1} below is the main result, regarding asymptotic expansions of  ``eigenvalue counting functions''. The other results below in the Introduction are  deduced from   Theorem \ref{main_thm1}.

\subsection*{Direct spectral problem}
Here,
we first introduce the following theorem, regarding asymptotic expansions of a kind of eigenvalue counting functions, where we use multi-index notations with
$$\alpha=(\alpha_1,\alpha_2,\dots,\alpha_{m-1})\in\left(\N\cup\{0\}\right)^{m-1},\quad \text{and}\,\,\, \beta=(1,2,\dots,{m-1}).$$  Also, we use $|\alpha|=\alpha_1+\alpha_2+\cdots+\alpha_{m-1}$, $\alpha!= \alpha_1!\alpha_2!\cdots\alpha_{m-1}!$ and $a^{\alpha}=a_1^{\alpha_1}a_2^{\alpha_2}\cdots a_{m-1}^{\alpha_{m-1}}$.
\begin{theorem}\label{main_thm1}
For  $a\in\C^{m-1}$, the eigenvalues $\lambda_n$ of $H_{\ell,P}$ satisfy
\begin{equation}\label{main_result}
\left(2n+1\right)\pi i=\sum_{j=0}^{\lfloor\frac{m+2}{2}\rfloor}d_{\ell,j}(a)\lambda_{n}^{\frac{1}{2}-\frac{j-1}{m}}+o(1),\quad\text{as $n\to\infty$},
\end{equation}
where $\lfloor x\rfloor$ is the largest integer that is less than or equal to $x$, and where
 the error $o(1)$ term is uniform on any compact set of $a\in\C^{m-1}$, and 
\begin{equation}\label{dmlj_def}
d_{\ell,j}(a)=\left\{
                    \begin{array}{cl}
                    2i\sqrt{\pi}\sin\left(\frac{\ell\pi}{m}\right)\frac{\Gamma\left(1+\frac{1}{m}\right)}{\Gamma\left(\frac{3}{2}+\frac{1}{m}\right)} \,\,\, &\text{if $j=0$,}\\
&\\
                   -4i\sum_{k=1}^j(-1)^{(\ell+1)k}K_{m,j,k}b_{j,k}(a)\sin\left(\frac{(j-1)\ell\pi}{m}\right)\cos\left(\frac{(j-1)\pi}{m}\right)\,\,\, &\text{if $1\leq j\leq\frac{m+1}{2}$,}\\
&\\
\eta_{m,\ell}(a)\,\,\, &\text{if $j=\frac{m+2}{2}$,}\\
                    \end{array}\right.
\end{equation}
where 
\begin{equation}\label{bjk_def}
b_{j,k}(a)={\frac{1}{2}\choose{k}}\sum_{\substack{|\alpha|=k\\ \alpha\cdot\beta=j}}\frac{k!}{\alpha!}\,a^{\alpha},\quad 1\leq k\leq j\leq \frac{m+2}{2},
\end{equation}
\begin{equation}\nonumber
\eta_{m,\ell}(a)=\left\{
                    \begin{array}{cl}
                    (-1)^{\frac{\ell-1}{2}}\frac{4\pi i}{m}\sum_{k=1}^{\frac{m+2}{2}}b_{\frac{m+2}{2},k}(a) \quad &\text{if $m$ is even and $\ell$ is odd,}\\
&\\
                   0\quad &\text{ otherwise,}
                    \end{array}\right.
\end{equation} 
and
\begin{equation}\nonumber
K_{m,j,k}
=\left\{
                    \begin{array}{cl}
-\frac{2}{m}
\quad &\text{if $j=k=1$},\\
&\\
                  -\frac{2k-1}{m+2-2j}B\left(k-\frac{j-1}{m},\,\frac{1}{2}+\frac{j-1}{m}\right)   \quad &\text{if $1\leq k\leq j\leq\frac{m+1}{2}$, $j\not=1$},\\
&\\
                 \frac{2}{m}\left(\ln 2-\frac{1}{1}-\frac{1}{3}-\dots-\frac{1}{2k-5}-\frac{1}{2k-3}\right)  \quad &\text{if $m$ is even, $1\leq k\leq j=\frac{m+2}{2}$,}
                    \end{array}\right.\nonumber
\end{equation}
where $B(\cdot,\cdot)$ is the beta function.
\end{theorem}
We obtain  \eqref{main_result} by investigating the asymptotic expansions of an entire function whose zeros are the eigenvalues.  Sibuya \cite{Sibuya} got \eqref{Bender_exp} by using the first order asymptotic expansion of the entire function. 

Next, we let $N_{\ell}(t)$, $t\in\R$, be the eigenvalue counting function, that is, $N_{\ell}(t)$ is the number of eigenvalues $\lambda$ of $H_{\ell,P}$ such that $|\lambda|\leq t$.  Then the following theorem on an asymptotic expansion of the eigenalue counting function is a consequence of Theorem \ref{main_thm1}.
\begin{theorem}\label{main_thm2}
Let $a\in\C^{m-1}$ be fixed. Suppose that $\Re\left(d_{\ell,j}(a)\right)=0$ for $1\leq j\leq\frac{m+2}{2}$. Then $N_{\ell}(t)$ has the asymptotic expansion 
\begin{equation}\label{count_ft}
N_{\ell}(t)=\frac{1}{2\pi i}\left(\sum_{j=0}^{\lfloor\frac{m+2}{2}\rfloor}d_{\ell,j}(a)t^{\frac{1}{2}-\frac{j-1}{m}}-\pi i\right)+O(1),\quad\text{as $t\to\infty$,}
\end{equation}
 where the error $O(1)$ is uniform for any compact set of $a\in\C^{m-1}$. 
\end{theorem}
\begin{proof}
In Theorem \ref{ineq_eq} below, we show that $|\lambda_n|<|\lambda_{n+1}|$ for all large $n\in\N$.

Suppose that $|\lambda_n|\leq t<|\lambda_{n+1}|$. Then since
$$\left(n+1+\frac{1}{2}\right)^{\frac{2m}{m+2}}=\left(n+\frac{1}{2}\right)^{\frac{2m}{m+2}}+O\left(n^{\frac{m-2}{m+2}}\right), \quad\text{as $n\to\infty$},$$
we see from Theorem \ref{eigen_asy} below that $|\lambda_{n+1}|-|\lambda_{n}|=O\left(n^{\frac{m-2}{m+2}}\right)$.
Thus, 
$$\lambda_n^{\frac{1}{2}-\frac{j-1}{m}}=t^{\frac{1}{2}-\frac{j-1}{m}}\left(1-\frac{t-\lambda_n}{t}\right)^{\frac{1}{2}-\frac{j-1}{m}}=t^{\frac{1}{2}-\frac{j-1}{m}}\left(1+O\left(\frac{t-\lambda_n}{t}\right)\right)=t^{\frac{1}{2}-\frac{j-1}{m}}+O\left(1\right).$$
Hence, replacing $\lambda_n^{\frac{1}{2}-\frac{j-1}{m}}$ in \eqref{main_result} by $t^{\frac{1}{2}-\frac{j-1}{m}}+O\left(1\right)$, and  solving the resulting equation for $n$ complete the proof.
\end{proof}

 Next, we improve the asymptotic expansion \eqref{Bender_exp} of the eigenvalues $\lambda_{n}$ of $H_{1,P}$. In particular, we will prove the following, which essentially invert \eqref{main_result} to get $\lambda_n$ in terms of $n$.
\begin{theorem}\label{eigen_asy}
For each $a\in\C^{m-1}$, there exist some constants $e_{j}(a)\in\C$, $2\leq j\leq \frac{m+2}{2}$, such that
\begin{equation}\label{main_asyeq}
\lambda_{n}=\lambda_{n,0}+\sum_{j=2}^{\lfloor\frac{m+2}{2}\rfloor}e_{j}(a)\lambda_{n,0}^{1-\frac{j}{m}}+o\left(\lambda_{n,0}^{1-\frac{1}{m}\lfloor\frac{m+2}{2}\rfloor}\right),\quad\text{as $n\to+\infty$},
\end{equation}
where the error term is uniform for any compact set of $a\in\C^{m-1}$, and where
\begin{equation}\nonumber
\lambda_{n,0}=\left(\frac{\sqrt{\pi}\Gamma\left(\frac{3}{2}+\frac{1}{m}\right)}{\sin\left(\frac{\ell\pi}{m}\right)\Gamma\left(1+\frac{1}{m}\right)}\left(n+\frac{1}{2}\right)\right)^{\frac{2m}{m+2}},
\end{equation}
and $e_j(a)$, $1\leq j\leq \frac{m+2}{2}$, are defined recurrently by $e_1(a)=0$ and
\begin{align}
&e_{j}(a)=-\frac{2m}{m+2}\left(\frac{d_{\ell,j}(a)}{d_{\ell,0}(a)}+\sum_{\substack{|\alpha|=k\geq 2\\ \alpha\cdot\beta=j}}{\frac{1}{2}+\frac{1}{m}\choose k}\frac{k!}{\alpha!}\, e(a)^{\alpha} +\sum_{r=2}^{j-2}\frac{d_{\ell,r}(a)}{d_{\ell,0}(a)}\sum_{\substack{|\alpha|=k\\ \alpha\cdot\beta=j-r}}{\frac{1}{2}+\frac{1-r}{m}\choose k}\frac{k!}{\alpha!}\, e(a)^{\alpha}\right),\nonumber
\end{align}
where $e(a)=(e_{1}(a),e_{2}(a),\dots,e_{m-1}(a))$.
\end{theorem}
We note for the first summation in the definition of $e_{j}(a)$ that $\alpha\cdot\beta=j$ implies $\alpha_i=0$ whenever $i\geq j$. Also, for the second summation, we point out that $\alpha\cdot\beta=j-r\leq j-2$ implies $\alpha_i=0$ whenever $i\geq j-1$.

The asymptotic expansions of the eigenvalues of $H_{\ell,P}$ with $\ell=\lfloor\frac{m}{2}\rfloor$ and $\ell=1$ have been studied by a number of people. For example,
 Maslov \cite{MAS}  computed the first three terms of asymptotic expansions of $\lambda_n^{\frac{3}{4}}$, where $\lambda_n$ are the eigenvalues of 
\begin{equation}\nonumber
-\frac{d^{2}}{dx^{2}}u+x^{4}u=\lambda u,\quad u\in L^2(\R). 
\end{equation}

Helffer and Robert \cite{HR} considered 
\begin{equation}\nonumber
-\frac{d^{2k}}{dx^{2k}}u+x^{2m}u+p(x)u=\lambda u,\quad u\in L^2(\R),
\end{equation}
where $k,\, m$ are positive integers and where $p(\cdot)$ is a {\it real}  polynomial of degree at most $2m-1$. They obtained existence of asymptotic expansions of eigenvalues to all orders, and suggested an explicit way of computing the coefficients of the asymptotic expansion. In particular, for the case when the potential is $\varepsilon x^4+x^2$, $\varepsilon>0$, Helffer and Robert \cite{HR} computed the first nine terms of the asymptotic expansion of $\lambda_n^{\frac{3}{4}}$. 

Fedoryuk \cite[\S 3.3]{Fedoryuk} considered \eqref{ptsym} with complex polynomial potentials and with \eqref{bdcond} for $\ell=\lfloor\frac{m}{2}\rfloor$, and showed the existence of asymptotic expansions of the eigenvalues to all orders. 
Note that there appear to be typographical errors in \cite[\S 3.3]{Fedoryuk}.  For example, when $m$ is even (and $\ell=\frac{m}{2}$) the leading coefficient of the asymptotic expansion of $\lambda_n$ in \cite[\S 3.3]{Fedoryuk}  is
$\left(2\frac{\sqrt{\pi}\Gamma\left(\frac{3}{2}+\frac{1}{m}\right)}{\Gamma\left(1+\frac{1}{m}\right)}\right)^{\frac{m+2}{2m}}$, which  is different from  $\left(\frac{\sqrt{\pi}\Gamma\left(\frac{3}{2}+\frac{1}{m}\right)}{\Gamma\left(1+\frac{1}{m}\right)}\right)^{\frac{2m}{m+2}}$ found in \cite{AA} and again in Theorem \ref{eigen_asy} above. 

Next, we point out some differences between work of Fedoryuk \cite[\S 3.3]{Fedoryuk} and the present work.  Fedoryuk \cite[\S 3.3]{Fedoryuk} showed the existence of asymptotic expansion of the eigenvalues to all orders while we do not. On the other hand, we treat all decaying boundary conditions \eqref{bdcond} with $1\leq\ell\leq m-1$,
while Fedoryuk \cite[\S 3.3]{Fedoryuk} studied the case  $\ell=\lfloor\frac{m}{2}\rfloor$ only. Moreover, we computed more coefficients $e_j(a)$ explicitly, and our methods  are different from Fedoryuk's.
  
\subsection*{Inverse spectral problem}
Here, we introduce results on inverse spectral problems, but first the following corollary is an easy consequence of Theorems \ref{main_thm1} and \ref{eigen_asy}, regarding how the coefficients of the asymptotic expansions depend on $a\in\C^{m-1}$.
\begin{corollary}\label{corollary_1}
Let $1\leq j\leq \frac{m+2}{2}$ be  a fixed integer. Then we have the following.
\begin{enumerate}
\item[(i)] $d_{\ell,j}(a)$ and $e_{j}(a)$ are  polynomials in $a_1,a_2,\dots,a_{j-1}, a_j$.
\item[(ii)]  $d_{\ell,j}(a)$ and $e_{j}(a)$ do not depend on $a_{j+1},a_{j+2},\dots, a_{m-1}$.
\item[(iii)] If $(j-1)\ell$ is a multiple of $m$, then $d_{\ell,j}(a)\equiv 0$, and $e_{j}(a)$ does not depend on $a_j$.
\item[(iv)] If $(j-1)\ell$ is not a multiple of $m$, then $d_{\ell,j}(a)$ and $e_{j}(a)$ depend linearly on $a_j$.
\end{enumerate}
\end{corollary}
\begin{proof}
Statements on $d_{\ell,j}(a)$ are direct consequences of the definition of $d_{\ell,j}(a)$ in Theorem \ref{main_thm1}. One can use statements on $d_{\ell,j}(a)$ and induction on $j$ to prove  statements on $e_{j}(a)$.
\end{proof}
Next, one can reconstruct some coefficients of the polynomial potential from the asymptotic expansion of the eigenvalues.
\begin{theorem}
 Let $1\leq j\leq\frac{m+1}{2}$ be a fixed integer. Suppose that $a_k$ is known whenever  $1\leq k\leq j$ and $(k-1)\ell$ is a multiple of $m$. If $(j-1)\ell$ is a multiple of $m$, then the asymptotic expansions of the eigenvalues $\lambda_{n}$ of $H_{\ell,P}$ of type \eqref{main_asyeq} with an error term $o\left(n^{\frac{2m-2j}{m+2}}\right)$ uniquely and explicitly determine $a_j$. 
\end{theorem}
\begin{proof}
From the asymptotic expansion of the eigenvalues, one gets $e_2(a),e_3(a),\dots,e_j(a)$ that are explicit polynomials in $a_1, a_2,\dots,a_j$.
Then since we know $a_k$ if $(k-1)\ell$ is a multiple of $m$, Corollary \ref{corollary_1} says that we can find  all $a_1, a_2,\dots,a_j$.
\end{proof}
\subsection*{Applications to $\mathcal{PT}$-symmetric potentials}
One says that $H_{\ell,P}$ is $\mathcal{PT}$-symmetric if the potential $V$ satisfies $\overline{V(-\overline{z})}=V(z)$, $z\in\C$, that is equivalent to $a\in\R^{m-1}$. Here, we prove the partial reality of the eigenvalues $\lambda_n$ of $\mathcal{PT}$-symmetric $H_{\ell,P}$. But first, we show the following theorem, regarding monotonicity of modulus of $\lambda_n$ for all large $n\in\N$.
\begin{theorem}\label{ineq_eq}
For each $a\in\C^{m-1}$ there exists $M>0$ such that $|\lambda_{n}|<|\lambda_{n+1}|$ if $n\geq M$. 
\end{theorem}
\begin{proof}
See  Theorem 3 in \cite{Shin2} for the proof of the case when $\ell=1$. One can see that proof of Theorem 3 in \cite{Shin2} can be easily adapted for the cases when $2\leq\ell\leq m-1$.
\end{proof}

Now we are ready to  prove the following theorem on the partial reality of the eigenvalues $\lambda_n$ of $H_{\ell,P}$.
\begin{theorem}\label{main_theorem}
Suppose that $a\in\R^{m-1}$. Then 
all but finitely many eigenvalues of $H_{\ell,P}$ are real and positive. Hence $\Re\left(d_{\ell,j}(a)\right)=0$ for all $1\leq j\leq\frac{m+2}{2}$, so that the counting function formula in Theorem \ref{main_thm2} is valid.
\end{theorem}
\begin{proof}
When $H_{\ell,P}$ is $\mathcal{PT}$-symmetric (i.e., $a\in\R^{m-1}$),  $u(z,\lambda)$ is an eigenfunction associated with an eigenvalue $\lambda$ if and only if $\overline{u(-\overline{z},\lambda)}$ is an eigenfunction associated with the eigenvalue $\overline{\lambda}$. Thus, the eigenvalues either appear in complex conjugate pairs, or else are real. So Theorem \ref{ineq_eq} implies Theorem \ref{main_theorem}.
\end{proof}
In recent years, these $\mathcal{PT}$-symmetric operators have gathered considerable attention, because ample numerical and asymptotic studies suggest that many of such operators have real eigenvalues only even though they are not self-adjoint. In particular, the differential operators $H$ with some polynomial potential $V$ and with the boundary condition \eqref{bdcond} have been considered in \cite{Bender,Bender2,CGM,KS,Ali1,Shin1,Shin,Simon,Znojil} and references therein. The rigorous proof of reality and positivity of the eigenvalues of  $\mathcal{PT}$-symmetric operators with certain classes of polynomial potentials and with the boundary condition \eqref{bdcond} for $\ell=1$,
was given by Dorey, Dunning and Tateo \cite{Dorey} in 2001 and by the present author \cite{Shin} in 2002. 

However, there are some $\mathcal{PT}$-symmetric polynomial potentials that produce a finite number of non-real eigenvalues \cite{Bender-1,Pham,Delabaere,Handy2,Handy1} for some particular classes of polynomial potentials.  So without any further restrictions on the real coefficients $a_k$, Theorem \ref{main_theorem} is the most general result one can expect about reality of eigenvalues of $\mathcal{PT}$-symmetric operators with polynomial potentials.

This paper is organized as follows. In Section \ref{prop_sect}, we will introduce work of Hille \cite{Hille} and Sibuya \cite{Sibuya}, regarding properties of solutions of \eqref{ptsym}. Also, we introduce  entire functions $W_{-1,\ell}(a,\lambda)$ whose zeros are closely related with the eigenvalues of $H_{\ell,P}$, due to Sibuya \cite{Sibuya} (c.\ f., Section \ref{asymp_eigen}). In Section \ref{sec_4}, we then provide asymptotics of the entire function $W_{-1,1}(a,\lambda)$ as $\lambda\to\infty$ in the complex plane \cite{Shin2}, improving the  asymptotics of $W_{-1,1}(a,\lambda)$ in \cite{Sibuya}.  In Section \ref{sec_5}, we  provide asymptotic expansions of $W_{-1,\ell}(a,\lambda)$ as $\lambda\to\infty$ in $\C$. In Section \ref{asymp_eigen}, we investigate how the zeros of $W_{-1,\ell}(a,\cdot)$ are related with the eigenvalues of $H_{\ell,P}$. In Sections \ref{sec_7} and \ref{sec_8}, we prove Theorem \ref{main_thm1}. In Section \ref{sec_9}, we prove Theorem \ref{eigen_asy}. Finally, in the Appendix we compute $K_{m,j,k}$ in Theorem \ref{main_thm1}, that is originally given in terms of certain integrals \eqref{K_definition}.

\section{Properties of the solutions}
\label{prop_sect}
In this section, we introduce work of Hille \cite{Hille} and Sibuya \cite{Sibuya} about properties of the solutions of \eqref{ptsym}.

First, we scale equation \eqref{ptsym} because many facts that we need later are stated for the scaled equation. Let $u$ be a solution of (\ref{ptsym}) and let $v(z,\lambda)=u(-iz,\lambda)$. 
 Then $v$ solves
\begin{equation}\label{rotated1}
-v^\dd(z,\lambda)+[(-1)^{\ell+1}z^m+P(z)+\lambda]v(z,\lambda)=0,
\end{equation}
where $m\geq 3$ and $P$  is a  polynomial (possibly, $P\equiv 0$) of the form
$$P(z)=a_1z^{m-1}+a_2z^{m-2}+\cdots +a_{m-1}z,\quad a_k\in\C.$$

When $\ell$ is odd, \eqref{rotated1} becomes 
\begin{equation}\label{rotated}
-v^\dd(z,\lambda)+[z^m+P(z)+\lambda]v(z,\lambda)=0.
\end{equation} 
Later we will treat the case when $\ell$ is even.

\begin{figure}[t]
    \begin{center}
    \includegraphics[width=.4\textwidth]{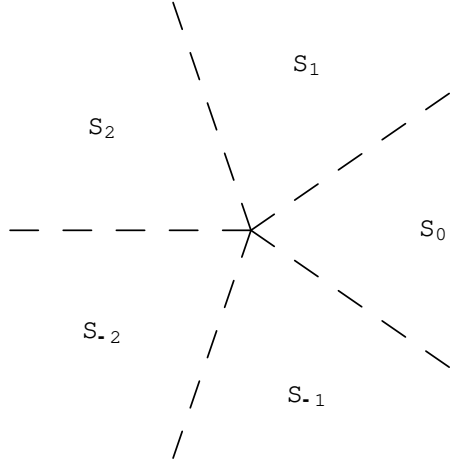}
    \end{center}
 \vspace{-.5cm}
\caption{The Stokes sectors for $m=3$. The dashed rays represent $\arg z=\pm\frac{\pi}{5},\,\pm\frac{3\pi}{5},\, \pi.$}\label{f:graph1}
\end{figure}

Since we scaled the argument of $u(\cdot,\lambda)$, we must rotate the boundary condition. We state them in a more general context by using the following definition.
\begin{definition}
{\rm {\it The Stokes sectors} $S_k$ of the equation (\ref{rotated})
are
$$ S_k=\left\{z\in \C:\left|\arg (z)-\frac{2k\pi}{m+2}\right|<\frac{\pi}{m+2}\right\}\quad\text{for}\quad k\in \Z.$$ }
\end{definition}
See Figure \ref{f:graph1}.

It is known from Hille \cite[\S 7.4]{Hille} that every nonconstant solution of (\ref{rotated}) either decays to zero or blows up exponentially, in each Stokes sector $S_k$. 
 That is, one has the following result.
\begin{lemma}[\protect{\cite[\S 7.4]{Hille}}]\label{gen_pro}
${} $

\begin{itemize}
\item [(i)] For each $k\in\Z$, every solution $v$ of (\ref{rotated})  is asymptotic to 
\begin{equation} \label{asymp-formula}
(const.)z^{-\frac{m}{4}}\exp\left[\pm \int^z \left[\xi^m+P(\xi)+\lambda\right]^{\frac{1}{2}}\,d\xi\right] 
\end{equation}
as $z \rightarrow \infty$ in every closed subsector of $S_k$.

\item [(ii)] If a nonconstant solution $v$ of \eqref{rotated} decays in $S_k$, it must blow up in $S_{k-1}\cup S_{k+1}$. However, when $v$ blows up in $S_k$, $v$ need not be decaying in $S_{k-1}$ or in $S_{k+1}$.
\end{itemize}
\end{lemma}
Lemma \ref{gen_pro} (i) implies that if $v$ decays along one ray in $S_k$, then it decays along all rays in $S_k$. Also, if $v$ blows up along one ray in $S_k$, then it blows up along all rays in $S_k$. This is essentially why we said in the Introduction that the boundary conditions \eqref{bdcond} with $1\leq\ell\leq m-1$ represent all decaying boundary conditions.

Still with $\ell$ odd, the two rays in \eqref{bdcond} map, by $z\mapsto -iz$,  to the rays $\arg(z)=\pm\frac{(\ell+1)\pi}{m+2}$ which are the center rays of the Stokes sectors $S_{\frac{\ell+1}{2}}$ and $S_{-\frac{\ell+1}{2}}$. Thus, the boundary conditions \eqref{bdcond} on $u$ become
\begin{equation}\nonumber
\text{$v(\cdot,\lambda)$ decays to zero in the Stokes sector $S_{\frac{\ell+1}{2}}$ and $S_{-\frac{\ell+1}{2}}$.}
\end{equation}

When $\ell$ is even, we let
 $y(z,\lambda)=v(\omega^{-\frac{1}{2}}z,\lambda)$ so that \eqref{rotated1} becomes
\begin{equation}\label{rotated2}
-y^\dd(z,\lambda)+[z^m+\omega^{-1} P(\omega^{-\frac{1}{2}}z)+\omega^{-1}\lambda]y(z,\lambda)=0,
\end{equation}
 where
$$\omega=\exp\left[\frac{2\pi i}{m+2}\right]$$
and hence, $\omega^{-\frac{m+2}{2}}=-1$.
For these cases, the boundary conditions \eqref{bdcond} become
\begin{equation}\nonumber
\text{$y(\cdot,\lambda)$ decays to zero in the Stokes sector $S_{\frac{\ell+2}{2}}$ and $S_{-\frac{\ell}{2}}$.}
\end{equation}

Next we will introduce Sibuya's results, but first we define  a sequence of complex numbers $b_j$ in terms of the $a_k$ and $\lambda$, as follows. For $\lambda \in \C$ fixed, we expand
\begin{eqnarray}
&&(1+a_1z^{-1}+a_2z^{-2}+\cdots+a_{m-1}z^{1-m}+\lambda z^{-m})^{1/2}\nonumber\\
&=&1+\sum_{k=1}^{\infty}{\frac{1}{2}\choose{k}}\left( a_1z^{-1}+a_2z^{-2}+\cdots+a_{m-1}z^{1-m}+\lambda z^{-m}\right)^k\nonumber\\
&=&1+\sum_{j=1}^{\infty}\frac{b_j(a,\lambda)}{z^j},\qquad\text{for large}\quad|z|.\label{b_def}
\end{eqnarray}
Note that $b_1,\,b_2,\,\ldots,\,b_{m-1}$ do not depend on $\lambda$, so we write $b_j(a)=b_j(a,\lambda)$ for $j=1,\,2,\dots,\,m-1$. So the above expansion without the $\lambda z^{-m}$ term still gives $b_j$ for $1\leq j\leq m-1$. We further define $r_m=-\frac{m}{4}$ if $m$ is odd, and $r_m=-\frac{m}{4}-b_{\frac{m}{2}+1}(a)$ if $m$ is even.

 The following theorem is a special case of Theorems 6.1, 7.2, 19.1 and 20.1 of Sibuya \cite{Sibuya} that is the main ingredient of the proofs of the main results in this paper. 
%The coefficient vector $a$ of $P$ is allowed to be complex, here. 
\begin{theorem}\label{prop}
Equation (\ref{rotated}), with $a\in \C^{m-1}$, admits a solution  $f(z,a,\lambda)$ with the following properties.
\begin{enumerate}
\item[(i)] $f(z,a,\lambda)$ is an entire function of $z,a $ and $\lambda$.
\item[(ii)] $f(z,a,\lambda)$ and $f^\d(z,a,\lambda)=\frac{\partial}{\partial z}f(z,a,\lambda)$ admit the following asymptotic expansions. Let $\varepsilon>0$. Then
\begin{align}
f(z,a,\lambda)=&\qquad z^{r_m}(1+O(z^{-1/2}))\exp\left[-F(z,a,\lambda)\right],\nonumber\\
f^\d(z,a,\lambda)=&-z^{r_m+\frac{m}{2}}(1+O(z^{-1/2}))\exp\left[-F(z,a,\lambda) \right],\nonumber
\end{align}
as $z$ tends to infinity in  the sector $|\arg z|\leq \frac{3\pi}{m+2}-\varepsilon$, uniformly on each compact set of $(a,\lambda)$-values . 
Here
\begin{equation}\nonumber
F(z,a,\lambda)=\frac{2}{m+2}z^{\frac{m}{2}+1}+\sum_{1\leq j<\frac{m}{2}+1}\frac{2}{m+2-2j}b_j(a) z^{\frac{1}{2}(m+2-2j)}.
\end{equation}
\item[(iii)] Properties \textup{(i)} and \textup{(ii)} uniquely determine the solution $f(z,a,\lambda)$ of (\ref{rotated}).
\item[(iv)] For each fixed $a\in\C^{m-1}$ and $\delta>0$, $f$ and $f^\d$ also admit the asymptotic expansions, 
\begin{align}
f(0,a,\lambda)=&[1+o(1)]\lambda^{-1/4}\exp\left[L(a,\lambda)\right],\label{eq1}\\
f^\d(0,a,\lambda)=&-[1+o(1)]\lambda^{1/4}\exp\left[L(a,\lambda)\right],\label{eq2}
\end{align}
as $\lambda\to\infty$ in the sector $|\arg(\lambda)|\leq\pi-\delta$, uniformly on each compact set of $a\in\C^{m-1}$, where 
\begin{align}
L(a,\lambda)=\left\{
                    \begin{array}{rl}
                    &\int_0^{+\infty}\left(\sqrt{t^m+P(t)+\lambda}- t^{\frac{m}{2}}-\sum_{j=1}^{\frac{m+1}{2}}b_j(a)t^{\frac{m}{2}-j}\right)\,dt \quad \text{if $m$ is odd,}\\
                   &\int_0^{+\infty}\left(\sqrt{t^m+P(t)+\lambda}- t^{\frac{m}{2}}-\sum_{j=1}^{\frac{m}{2}}b_j(a)t^{\frac{m}{2}-j}-\frac{b_{\frac{m}{2}+1}}{t+1}\right)\,dt  \quad \text{if $m$ is even.}
                    \end{array}
                         \right. \nonumber
\end{align}
\item[(v)] The entire functions  $\lambda\mapsto f(0,a,\lambda)$ and $\lambda\mapsto f^\d(0,a,\lambda)$ have orders $\frac{1}{2}+\frac{1}{m}$.
\end{enumerate}
\end{theorem}
\begin{proof}
In  Sibuya's book \cite{Sibuya}, see Theorem 6.1 for a proof of (i) and (ii); Theorem 7.2 for a proof of (iii); and Theorem 19.1 for a proof of (iv).  Moreover, (v) is a consequence of (iv) along with Theorem 20.1 in \cite{Sibuya}. Note that properties (i), (ii) and (iv) are summarized on pages 112--113 of Sibuya \cite{Sibuya}. 
\end{proof}

\begin{remarks}
{\rm (I) Uniformness of the error term in Theorem \ref{main_thm1}  is essentially due to uniformness of error terms in \eqref{eq1} and \eqref{eq2}. One can check this by carefully following our proofs. In this paper, we omit this part of the proof.

(II) Throughout this paper, we will deal with numbers like $\left(\omega^{\nu}\lambda\right)^{s}$ for some $s\in\R$, and $\nu\in\C$. As usual, we will use
$$\omega^{\nu}=\exp\left[\nu \frac{2\pi i}{m+2}\right]$$
and if $\arg(\lambda)$ is specified, then
$$\arg\left(\left(\omega^{\nu}\lambda\right)^{s}\right)=s\left[\arg(\omega^{\nu})+\arg(\lambda)\right]=s\left[\Re(\nu)\frac{2\pi}{m+2}+\arg(\lambda)\right],\quad s\in\R.$$
If $s\not\in\Z$ then the branch of $\lambda^s$ is chosen to be the negative real axis.
}
\end{remarks}

Next, we provide an asymptotic expansion of $L$ in \cite{Shin2}.   But first, we recall that 
for $1\leq k\leq j\leq \frac{m+2}{2}$,
\begin{equation}\nonumber
b_{j,k}(a)={\frac{1}{2}\choose{k}}\sum_{\substack{j_p\not=j_q\text{if}\, p\not=q,\, i_p\geq 1\\ i_1+\cdots+i_s=k\\ i_1j_1+\cdots+i_sj_s=j}}\frac{k!}{i_1!i_2!\cdots i_s!}a_{j_1}^{i_1}a_{j_2}^{i_2}\cdots a_{j_{s}}^{i_{s}}={\frac{1}{2}\choose{k}}\sum_{\substack{|\alpha|=k\\ \alpha\cdot\beta=j}}\frac{k!}{\alpha!}\,a^{\alpha}.
\end{equation}
 Then
\begin{equation}\nonumber
b_j(a)=\sum_{k=1}^j b_{j,k}(a).
\end{equation}
Also, we define
\begin{equation}\label{K_definition}
K_{m,j,k}=\left\{\begin{array}{rl}
&\int_0^{\infty}\left(\frac{t^{mk-j}}{\left(t^m+1\right)^{k-\frac{1}{2}}}-t^{\frac{m}{2}-j}\right)\,dt,\,\,1\leq k\leq j\leq \frac{m}{2}\,\,\text{or}\,\,k=j=0,\\
&\int_0^{\infty}\left(\frac{t^{mk-\frac{m}{2}-1}}{\left(t^m+1\right)^{k-\frac{1}{2}}}-\frac{1}{t+1}\right)\,dt,\,\,1\leq k\leq j= \frac{m+2}{2} \,\,\text{for $m$ even,}
\end{array}
\right.
\end{equation}
and define
\begin{equation}\label{K_deF}
K_{m,0}(a)=K_{m}=K_{m,0,0},\quad K_{m,j}(a)=\sum_{k=1}^jb_{j,k}(a)K_{m,j,k},\quad 1\leq j\leq \frac{m+2}{2}.
\end{equation}
See Appendix for $K_{m,j,k}$ in terms of some gamma functions as in Theorem \ref{main_thm1}.
\begin{lemma}\label{asy_lemma}
Let $m\geq 3$ and $a\in\C^{m-1}$ be fixed. Then there exist constants $K_{m,j}(a)\in\C$, $0\leq j \leq\frac{m}{2}+1$, such that 
\begin{equation}
L(a,\lambda)=\left\{\begin{array}{rl}
&\sum_{j=0}^{\frac{m+1}{2}}K_{m,j}(a)\lambda^{\frac{1}{2}+\frac{1-j}{m}}+O\left(|\lambda|^{-\frac{1}{2m}}\right)\,\,\text{if $m$ is odd,}\\
&\sum_{j=0}^{\frac{m}{2}+1}K_{m,j}(a)\lambda^{\frac{1}{2}+\frac{1-j}{m}}-\frac{b_{\frac{m}{2}+1}(a)}{m}\ln(\lambda)+O\left(|\lambda|^{-\frac{1}{m}}\right)\,\,\text{if $m$ is even,}
\end{array}
\right.\nonumber
\end{equation}
as $\lambda\to\infty$ in the sector $|\arg(\lambda)|\leq \pi-\delta$, uniformly on each compact set of $a\in\C^{m-1}$. 
\end{lemma}
\begin{proof}
See \cite{Shin2} for a proof.
\end{proof} 
Sibuya \cite{Sibuya} proved the following corollary, directly from Theorem \ref{prop}, that will be used later in Sections \ref{sec_4} and \ref{sec_5}.
\begin{corollary}\label{lemma_decay}
Let $a\in\C^{m-1}$ be fixed. Then
$L(a,\lambda)=K_{m}\lambda^{\frac{1}{2}+\frac{1}{m}}(1+o(1))$ as $\lambda$ tends to infinity in the sector $|\arg \lambda|\leq \pi-\delta$, and hence
\begin{equation}\label{re_part}
\Re \left(L(a,\lambda)\right)=K_{m}\cos\left(\frac{m+2}{2m}\arg(\lambda)\right)|\lambda|^{\frac{1}{2}+\frac{1}{m}}(1+o(1))
\end{equation}
as $\lambda\to\infty$ in the sector $|\arg (\lambda)|\leq \pi-\delta$.

In particular, $\Re \left(L(a,\lambda)\right)\to+\infty$ as $\lambda\to\infty$ in any closed subsector of the sector $|\arg(\lambda)|<\frac{m\pi}{m+2}$. In addition, $\Re \left(L(a,\lambda)\right)\to-\infty$ as $\lambda\to\infty$ in any closed subsector of the sectors $\frac{m\pi}{m+2}<|\arg(\lambda)|<\pi-\delta$.
\end{corollary}

Based on Corollary \ref{lemma_decay},
Sibuya \cite[Theorem 29.1]{Sibuya} also proved the  asymptotic expansion \eqref{Bender_exp} of the eigenvalues of $H_{\ell,P}$ for $\ell=1$.

Also, Sibuya \cite{Sibuya} constructed solutions of \eqref{rotated} that decays in $S_k$, $k\in\Z$. Before we introduce this, we let 
\begin{equation}\label{G_def}
G^{\ell}(a):=(\omega^{(m+1)\ell}a_1, \omega^{m\ell}a_2,\ldots,\omega^{3\ell}a_{m-1})\quad \text{for}\quad \ell\in \frac{1}{2}\Z.
\end{equation}
Then we have the following lemma, regarding properties of $G^{\ell}(\cdot)$.
\begin{lemma}\label{lemma_25}
For $a\in\C^{m-1}$ fixed, and $\ell_1,\ell_2,\ell\in\frac{1}{2}\Z$,
$G^{\ell_1}(G^{\ell_2}(a))=G^{\ell_1+\ell_2}(a)$, and
\begin{equation}\nonumber
b_{j,k}(G^{\ell}(a))=\omega^{((m+2)k-j)\ell}b_{j,k}(a),\quad \ell\in\frac{1}{2}\Z.
\end{equation}

If $\ell\in\Z$ then
\begin{equation}\nonumber
b_j(G^{\ell}(a))=\omega^{-j\ell}b_j(a).
\end{equation}
\end{lemma}

%%%%%%%%%%%%%%%%%%%%%%%%%%%%%%%%%%%%%%%%%%%%%%%%%%%
%%%%%%%%%%%%%%%%%%%%%%%%%%%%%%%%%%%%%%%%%%%%%%%%%%
%\section{Eigenvalues are zeros of an entire function}\label{entire_sect}
%%%%%%%%%%%%%%%%%%%%%%%%%%%%%%%%%%%%%%%%%%%%%%%%%%%%
%%%%%%%%%%%%%%%%%%%%%%%%%%%%%%%%%%%%%%%%%%%%%%%%%%%%%%%

%In this section, we will prove that the eigenvalues are 
%zeros of an entire function.
  
Next, recall that the function
 $f(z,a,\lambda)$  in Theorem \ref{prop} solves (\ref{rotated}) and decays to zero exponentially as $z\rightarrow \infty$ in $S_0$, and  blows up in $S_{-1}\cup S_1$. One can check that the function
$$f_k(z,a,\lambda):=f(\omega^{-k}z,G^k(a),\omega^{2k}\lambda),\quad k\in\Z,$$
 which is obtained by scaling $f(z,G^k(a),\omega^{2k}\lambda)$ in the $z$-variable, also solves (\ref{rotated}). It is clear that $f_0(z,a,\lambda)=f(z,a,\lambda)$, and that $f_k(z,a,\lambda)$ decays in $S_k$ and blows up in $S_{k-1}\cup S_{k+1}$ since $f(z,G^k(a),\omega^{2k}\lambda)$ decays in $S_0$. Since no nonconstant solution decays in two consecutive Stokes sectors (see Lemma \ref{gen_pro} (ii)), $f_{k}$ and $f_{k+1}$ are linearly independent and hence any solution of (\ref{rotated}) can be expressed as a linear combination of these two. Especially,  for each $k\in\Z$ there exist some coefficients $C_k(a,\lambda)$ and $\widetilde{C}_k(a,\lambda)$ such that
\begin{equation}\label{stokes}
f_{k}(z,a,\lambda)=C_k(a,\lambda)f_{0}(z,a,\lambda)+\widetilde{C}_k(a,\lambda)f_{-1}(z,a,\lambda).
\end{equation}
We then see that 
\begin{equation}\label{C_def}
C_k(a,\lambda)=-\frac{W_{k,-1}(a,\lambda)}{W_{-1,0}(a,\lambda)}\quad\text{and}\quad \widetilde{C}_k(a,\lambda)=\frac{W_{k,0}(a,\lambda)}{W_{-1,0}(a,\lambda)},
\end{equation}
where $W_{j,\ell}=f_jf_{\ell}^\d -f_j^\d f_{\ell}$ is the Wronskian of $f_j$ and $f_{\ell}$. Since both $f_j,\,f_{\ell}$ are solutions of the same linear equation (\ref{rotated}), we know that the Wronskians are constant functions of $z$. Also, $f_k$ and $f_{k+1}$ are linearly independent, and hence $W_{k,k+1}\not=0$ for all $k\in \Z$. 

Also, the following is an easy consequence of \eqref{stokes} and \eqref{C_def}. For each $k,\ell\in\Z$ we have 
\begin{align}
W_{\ell,k}(a,\lambda)&=C_k(a,\lambda)W_{\ell,0}(a,\lambda)+\widetilde{C}_k(a,\lambda)W_{\ell,-1}(a,\lambda)\nonumber\\
&=-\frac{W_{k,-1}(a,\lambda)W_{\ell,0}(a,\lambda)}{W_{-1,0}(a,\lambda)}+\frac{W_{k,0}(a,\lambda)W_{\ell,-1}(a,\lambda)}{W_{-1,0}(a,\lambda)}.\label{kplus}
\end{align}

Moreover, we have the following lemma that is useful  later on.
\begin{lemma}\label{shift_lemma}
Suppose $k,\,j\in\Z$. Then
\begin{equation}\label{kplus1}
W_{k+1,j+1}(a,\lambda)=\omega^{-1}W_{k,j}(G(a),\omega^2\lambda),
\end{equation}
and $W_{0,1}(a,\lambda)=2\omega^{\mu(a)}$, where
\begin{eqnarray}
\mu(a)=\left\{
              \begin{array}{rl}
              \frac{m}{4}  \quad &\text{if $m$ is odd,}\\
              \frac{m}{4}- b_{\frac{m}{2}+1}(a) \quad &\text{if $m$ is even.}
              \end{array}
                         \right. \nonumber
\end{eqnarray}

\end{lemma}
\begin{proof}
See Sibuya \cite[pages 116-118]{Sibuya}.
\end{proof}

We let $\nu(a)=\frac{m}{4}-\mu(a)$, that is, 
\begin{eqnarray}
\nu(a)=\left\{
                         \begin{array}{rl}
                         0 & \quad \text{if $m$ is odd,}\\
                         b_{\frac{m}{2}+1}(a)& \quad \text{if $m$ is even.}
                          \end{array}
                         \right. \label{def_nu}
\end{eqnarray}

%%%%%%%%%%%%%%%%%%%%%%%%%%%%%%%%%%%%%%%%%%%%%%
\section{Asymptotics of $W_{-1,1}(a,\lambda)$}\label{sec_4}
%%%%%%%%%%%%%%%%%%%%%%%%%%%%%%%%%%%%%%%%%%%%%%%%%%%
In this section, we introduce asymptotic expansions of $W_{-1,1}(a,\lambda)$ as $\lambda\to\infty$ along the rays in the complex plane \cite{Shin2}.

First, we provide an asymptotic expansion of  the Wronskian  $W_{0,j}(a,\lambda)$ of $f_0$ and $f_j$ that will be frequently used later.
\begin{lemma}\label{lemma7}
Suppose that $1\leq j\leq \frac{m}{2}+1$. Then for each $a\in\C^{m-1}$,
\begin{equation}\label{sec_eq1}
W_{0,j}(a,\lambda)=[2i\omega^{-\frac{j}{2}}+o(1)]\exp\left[L(G^{j}(a),\omega^{2j-m-2}\lambda)+L(a,\lambda)\right],
\end{equation}
as $\lambda\to\infty$ in the sector
\begin{equation}\label{sector0}
-\pi+\delta\leq \pi-\frac{4j\pi}{m+2}+\delta \leq \arg(\lambda)\leq \pi-\delta.
\end{equation}
\end{lemma}

Next, we provide an asymptotic expansion of $W_{-1,1}(a,\lambda)$ as $\lambda\to\infty$ in the sector near the negative real axis.
%%%%%%%%%%%%%%%%%%%%%%%%%%%%%%%%
\begin{theorem}\label{thm_neg}
Let $m\geq 3$, $a\in\C^{m-1}$ and $0<\delta<\frac{\pi}{m+2}$ be fixed. Then
\begin{equation}\label{asy_1}
W_{-1,1}(a,\lambda)=[2i+o(1)]\exp\left[L(G^{-1}(a),\omega^{-2}\lambda)+L(G(a),\omega^{-m}\lambda)\right],
\end{equation}
as $\lambda\to \infty$ along the rays in the sector
\begin{equation}\label{sector1}
\pi-\frac{4\pi}{m+2}+\delta\leq \arg(\lambda)\leq \pi+\frac{4\pi}{m+2}-\delta.
\end{equation}
\end{theorem}
%%%%%%%%%%%%%%%%%%%%%%%%%%%%%%%%%
\begin{proof}
This is an easy consequence of Lemma \ref{lemma7} with $j=2$ and \eqref{kplus1}.
\end{proof}

Also,  for integers $m\geq 4$ we provide an asymptotic expansion of $W_{-1,1}(a,\lambda)$ as $\lambda\to\infty$ in the sector $|\arg(\lambda)|\leq \pi-\delta$.
%%%%%%%%%%%%%%%%%%%%%%%%%%%%%%%%%%%%%%%%%%%%%%%%%%%%%%
\begin{theorem}\label{zero_thm}
Let $a\in\C^{m-1}$ and $0<\delta<\frac{\pi}{2(m+2)}$ be fixed. 
If $m\geq 4$ then 
\begin{align}
W_{-1,1}(a,\lambda)=&[2\omega^{\frac{1}{2}+\mu(a)}+o(1)]\exp\left[L(G^{-1}(a),\omega^{-2}\lambda)-L(a,\lambda)\right]\nonumber\\
&+[2\omega^{\frac{1}{2}+\mu(a)+2\nu(a)}+o(1)]\exp\left[L(G(a),\omega^{2}\lambda)-L(a,\lambda)\right],\label{tot_asy}
\end{align}
as $\lambda\to\infty$ in the sector 
\begin{equation}\label{sector41}
-\pi+\delta \leq \arg(\lambda)\leq \pi-\delta.
\end{equation} 
\end{theorem}
%%%%%%%%%%%%%%%%%%%%%%%%%%%%%%%%%%%%%%%%%%%%%

We provide an asymptotic expansion of $W_{-1,1}(a,\lambda)$ as $\lambda\to\infty$ along the rays in the upper half plane. 
%%%%%%%%%%%%%%%%%%%%%%%%%%%%%%%%%%%%%%%%%%%%
\begin{corollary}\label{lemma_up}
Let $m\geq 4$, $a\in\C^{m-1}$ and $0<\delta<\frac{\pi}{m+2}$ be fixed. Then
\begin{equation}\nonumber
W_{-1,1}(a,\lambda)=[2\omega^{\frac{1}{2}+\mu(a)}+o(1)]\exp\left[L(G^{-1}(a),\omega^{-2}\lambda)-L(a,\lambda)\right],
\end{equation}
as $\lambda\to\infty$ in the sector $\delta \leq \arg(\lambda)\leq \pi-\delta$.
Also,
\begin{equation}\nonumber
W_{-1,1}(a,\lambda)=[2\omega^{\frac{1}{2}+\mu(a)+2\nu(a)}+o(1)]\exp\left[L(G(a),\omega^{2}\lambda)-L(a,\lambda)\right],
\end{equation}
as $\lambda\to\infty$ in the sector $-\pi+\delta \leq \arg(\lambda)\leq -\delta$.
\end{corollary}
\begin{proof}
We will determine which term in \eqref{tot_asy} dominates in the upper and lower half planes.

Since, by \eqref{re_part},
\begin{equation}\nonumber
\Re(L(a,\lambda))=K_m\cos\left(\frac{m+2}{2m}\arg(\lambda)\right)|\lambda|^{\frac{1}{2}+\frac{1}{m}}(1+o(1)),
\end{equation}
we have
\begin{align}
&\left[\Re(L(G^{-1}(a),\omega^{-2}\lambda))-\Re(L(a,\lambda))\right]-\left[\Re(L(G(a),\omega^{2}\lambda))-\Re(L(a,\lambda))\right]\nonumber\\
&=K_m\left[\cos\left(-\frac{2\pi}{m}+\frac{m+2}{2m}\arg(\lambda)\right)-\cos\left(\frac{2\pi}{m}+\frac{m+2}{2m}\arg(\lambda)\right)\right]|\lambda|^{\frac{1}{2}+\frac{1}{m}}(1+o(1))\nonumber\\
&=2K_m\sin\left(\frac{2\pi}{m}\right)\sin\left(\frac{m+2}{2m}\arg(\lambda)\right)|\lambda|^{\frac{1}{2}+\frac{1}{m}}(1+o(1)).\nonumber
\end{align}
Thus, the first term in \eqref{tot_asy} dominates as $\lambda\to\infty$ along the rays in the upper half plane, and the second term dominates in the lower half plane. This completes the proof. 
\end{proof}
%%%%%%%%%%%%%%%%%%%%%%%%%%%%%%%%%%%%%%%%%%%%%%%%%%%%%%%%%%%%%%%%%%%%%%
%%%%%%%%%%%%%%%%%%%%%%%%%%%%%%%%%%%%%%%%%%%%%%%%%%%%%%%%%%%%%%%%%%%%%%
\begin{proof}[Proof of Theorem ~\ref{zero_thm}]
In \cite{Shin2}, $C(a,\lambda)$ is used for $\frac{W_{-1,1}(a,\lambda)}{W_{0,1}(a,\lambda)}$ and asymptotics of $C(a,\lambda)$ are provided.  Notice that $W_{-1,1}(a,\lambda)=2\omega^{\mu(a)}C(a,\lambda)$.

Theorem 13 in \cite{Shin2} implies
 \eqref{tot_asy}  for the sector 
\begin{equation}\label{sector4}
\pi-\frac{4\lfloor\frac{m}{2}\rfloor\pi}{m+2}+\delta \leq \arg(\lambda)\leq \pi-\frac{4\pi}{m+2}-\delta.
\end{equation} 
Theorem 14 in \cite{Shin2} implies that
\begin{align}
W_{-1,1}(a,\lambda)=&[2\omega^{\frac{1}{2}+\mu(a)}+o(1)]\exp\left[L(G^{-1}(a),\omega^{-2}\lambda)-L(a,\lambda)\right]\nonumber\\
&+[2\omega^{1+2\mu(a)+4\nu(a)}+o(1)]\exp\left[-L(G^2(a),\omega^{2-m}\lambda)-L(a,\lambda)\right],\nonumber
\end{align}
as $\lambda\to\infty$ in the sector $\pi-\frac{8\pi}{m+2}+\delta\leq\arg(\lambda)\leq\pi-\delta$. One can check that the first term dominates in this sector, by using an argument  similar to that in the proof of Corollary \ref{lemma_up}.

Also, Theorem 15 in \cite{Shin2} implies that
\begin{align}
W_{-1,1}(a,\lambda)=&[2\omega^{1+2\mu(a)}+o(1)]\exp\left[-L(a,\omega^{-m-2}\lambda)-L(^{-2}(a),\omega^{-4}\lambda)\right]\nonumber\\
&+[2\omega^{\frac{1}{2}+\mu(a)+2\nu(a)}+o(1)]\exp\left[L(G(a),\omega^{-m}\lambda)-L(a,\omega^{-m-2}\lambda)\right],\nonumber
\end{align}
as $\lambda\to\infty$ in the sector $\pi+\delta\leq\arg(\lambda)\leq\pi+\frac{8\pi}{m+2}-\delta$. One can check that the second term dominates in this sector. Then we replace $\lambda$ by $\omega^{m+2}\lambda$ to convert the sector here to $-\pi+\delta\leq\arg(\lambda)\leq-\pi+\frac{8\pi}{m+2}-\delta$. This completes the proof.
\end{proof}

\begin{theorem}\label{thm_sector2}
Let $m=3$ and let $a\in\C^{m-1}$ and $0<\delta<\frac{\pi}{m+2}$ be fixed. Then
\begin{align}
W_{-1,1}(a,\lambda)=&[-2\omega^{-\frac{5}{4}}+o(1)]\exp\left[L(G^{4}(a),\omega^{-2}\lambda)-L(a,\lambda)\right]\nonumber\\
&-[2i\omega^{\frac{5}{2}}+o(1)]\exp\left[-L(G^2(a),\omega^{-1}\lambda)-L(a,\lambda)\right],\nonumber
\end{align}
as $\lambda\to\infty$ in the sector 
$-\delta \leq \arg(\lambda)\leq \pi-\delta$.
Also,
\begin{align}
W_{-1,1}(a,\lambda)=&[-2i\omega^{\frac{5}{2}}+o(1)]\exp\left[-L(a,\omega^{-5}\lambda)-L(G^{-2}(a),\omega^{-4}\lambda)\right]\nonumber\\
&+[2\omega^{\frac{15}{4}}+o(1)]\exp\left[L(G(a),\omega^{-3}\lambda)-L(a,\omega^{-5}\lambda)\right],\nonumber
\end{align}
as $\lambda\to\infty$ in the sector $\pi+\delta \leq \arg(\lambda)\leq \delta$.
\end{theorem}
\begin{proof}
See Theorems 14 and 15 in \cite{Shin2} for a proof.
\end{proof}

%%%%%%%%%%%%%%%%%%%%%%%%%%%%%%%%%%%%%%%%%%%%%%%%
%%%%%%%%%%%%%%%%%%%%%%%%%%%%%%%%%%%%%%%%%%%%%%%%
\section{Asymptotics of $W_{-1,n}(a,\lambda)$}\label{sec_5}
%%%%%%%%%%%%%%%%%%%%%%%%%%%%%%%%%%%%%%%%%%%%%%%%
In this section, we will provide asymptotic expansions of $W_{-1,n}(a,\cdot)$, zeros of which will be closely related with the eigenvalues of $H_{n,P}$.

First, we treat the cases when  $1\leq n <\lfloor\frac{m}{2}\rfloor$.
\begin{theorem}
Let $1\leq n <\lfloor\frac{m}{2}\rfloor$ be an integer. Then $W_{-1,n}(a,\cdot)$  admits the following asymptotic expansion 
\begin{equation}\label{up_asy}
W_{-1,n}(a,\lambda)=[2\omega^{\frac{2-n}{2}+\mu(G^{n-1}(a))}+o(1)]\exp\left[L(G^{-1}(a),\omega^{-2}\lambda)-L(G^{n-1}(a),\omega^{2(n-1)}\lambda)\right],
\end{equation}
as $\lambda\to\infty$ in the sector 
\begin{equation}\label{zero3_sector}
-\frac{2(n-1)\pi}{m+2}+\delta\leq\arg(\lambda)\leq\pi -\frac{4n\pi}{m+2}+\delta.
\end{equation}

Also,
\begin{align}
W_{-1,n}(a,\lambda)&=[2\omega^{\frac{2-n}{2}+\mu(G^{n-1}(a))}+o(1)]\exp\left[L(G^{-1}(a),\omega^{-2}\lambda)-L(G^{n-1}(a),\omega^{2(n-1)}\lambda)\right]\nonumber\\
&+[2\omega^{\frac{2-n}{2}+\mu(G^{-1}(a))}+o(1)]\exp\left[L(G^n(a),\omega^{2n}\lambda)-L(a,\lambda)\right],\label{zero_asy}
\end{align}
as $\lambda\to\infty$ in the sector 
\begin{equation}\label{zero2_sector}
-\frac{2(n-1)\pi}{m+2}-\delta\leq\arg(\lambda)\leq-\frac{2(n-1)\pi}{m+2}+\delta.
\end{equation}
\end{theorem}
\begin{proof}
First we will prove \eqref{up_asy} for the sector
 \begin{equation}\label{zero1_sector}
-\frac{2(n-1)\pi}{m+2}+\delta\leq\arg(\lambda)\leq\pi -\frac{4n\pi}{m+2}-\delta
\end{equation}
 and the second part of the theorem by induction on $n$.

The case when $n=1$ is trivially satisfied by  Theorem \ref{zero_thm} and Corollary \ref{lemma_up} since $\mu(a)+2\nu(a)=\mu(G^{-1}(a))$.

Suppose that \eqref{up_asy} holds in the sector \eqref{zero1_sector} for $n-1$. From this induction hypothesis we have
\begin{align}
W_{0,n}(a,\lambda)&=\omega^{-1}W_{-1,n-1}(G(a),\omega^2\lambda)\nonumber\\
&=[2\omega^{-\frac{n-1}{2}+\mu(G^{n-1}(a))}+o(1)]\exp\left[L(a,\lambda)-L(G^{n-1}(a),\omega^{2(n-1)}\lambda)\right],\label{1steq}
\end{align}
as $\lambda\to\infty$ in the sector
\begin{equation}\nonumber
-\frac{2(n-2)\pi}{m+2}+\delta\leq\arg(\omega^2\lambda)\leq\pi -\frac{4(n-1)\pi}{m+2}-\delta,
\end{equation}
that is,
\begin{equation}\label{dom_sector}
-\frac{2n\pi}{m+2}+\delta\leq\arg(\lambda)\leq\pi -\frac{4n\pi}{m+2}-\delta.
\end{equation}

Also, from Lemma \ref{lemma7} if $1\leq j\leq\frac{m}{2}+1$, then we have 
\begin{equation}\label{2ndeq}
W_{0,j}(a,\lambda)=[2i\omega^{-\frac{j}{2}}+o(1)]\exp\left[L(G^{j}(a),\omega^{2j-m-2}\lambda)+L(a,\lambda)\right],
\end{equation}
as $\lambda\to\infty$ in the sector
\begin{equation}\nonumber
 \pi-\frac{4j\pi}{m+2}+\delta \leq \arg(\lambda)\leq \pi-\delta.
\end{equation}
We solve \eqref{kplus} for $W_{\ell,-1}(a,\lambda)$ and set $\ell=n$ to get
\begin{equation}\label{kn_asy}
W_{-1,n}(a,\lambda)=\frac{W_{-1,0}(a,\lambda)W_{n,k}(a,\lambda)}{W_{0,k}(a,\lambda)}+\frac{W_{0,n}(a,\lambda)W_{-1,k}(a,\lambda)}{W_{0,k}(a,\lambda)}
\end{equation}
Set $k=\lfloor\frac{m}{2}\rfloor$. Then since $1\leq k-n<k=\lfloor\frac{m}{2}\rfloor$, using \eqref{kplus1},
\begin{align}
W_{-1,n}(a,\lambda)&=\frac{2\omega^{\mu(G^{-1}(a))}W_{0,k-n}(G^n(a),\omega^{2n}\lambda)}{\omega^{n-1}W_{0,k}(a,\lambda)}+\frac{W_{0,n}(a,\lambda)W_{0,k+1}(G^{-1}(a),\omega^{-2}\lambda)}{\omega^{-1}W_{0,k}(a,\lambda)}\nonumber\\
&=\frac{2\omega^{\mu(G^{-1}(a))}[2i\omega^{-\frac{k-n}{2}}+o(1)]\exp\left[L(G^{n}(a),\omega^{2k-m-2}\lambda)+L(G^n(a),\omega^{2n}\lambda)\right]}{\omega^{n-1}[2i\omega^{-\frac{k}{2}}+o(1)]\exp\left[L(G^{k}(a),\omega^{2k-m-2}\lambda)+L(a,\lambda)\right]}\nonumber\\
&+\frac{[2\omega^{-\frac{n-1}{2}+\mu(G^{n-1}(a))}+o(1)]\exp\left[L(a,\lambda)-L(G^{n-1}(a),\omega^{2(n-1)}\lambda)\right]}{\omega^{-1}[2i\omega^{-\frac{k}{2}}+o(1)]\exp\left[L(G^{k}(a),\omega^{2k-m-2}\lambda)+L(a,\lambda)\right]}\nonumber\\
&\times [2i\omega^{-\frac{k+1}{2}}+o(1)]\exp\left[L(G^{k}(a),\omega^{2k-m-2}\lambda)+L(G^{-1}(a),\omega^{-2}\lambda)\right]\nonumber\\
&=[2\omega^{\frac{2-n}{2}+\mu(G^{-1}(a))}+o(1)]\exp\left[L(G^n(a),\omega^{2n}\lambda)-L(a,\lambda)\right]\nonumber\\
&+[2\omega^{\frac{2-n}{2}+\mu(G^{n-1}(a))}+o(1)]\exp\left[L(G^{-1}(a),\omega^{-2}\lambda)-L(G^{n-1}(a),\omega^{2(n-1)}\lambda)\right],\label{asy_sector}
\end{align}
where we used \eqref{1steq} for $W_{0,n}(a,\lambda)$ and \eqref{2ndeq} for everything else, provided that $\lambda$ lies in \eqref{dom_sector} and that
\begin{align}
\pi-\frac{4\left(\lfloor\frac{m}{2}\rfloor-n\right)\pi}{m+2}+\delta \leq &\arg(\omega^{2n}\lambda)\leq \pi-\delta\nonumber\\
\pi-\frac{4\lfloor\frac{m}{2}\rfloor\pi}{m+2}+\delta \leq &\arg(\lambda)\leq \pi-\delta\nonumber\\
\pi-\frac{4\left(\lfloor\frac{m}{2}\rfloor+1\right)\pi}{m+2}+\delta \leq &\arg(\omega^{-2}\lambda)\leq \pi-\delta,\nonumber
\end{align}
that is,
\begin{equation}\nonumber
-\frac{2n\pi}{m+2}+\delta\leq\arg(\lambda)\leq\pi -\frac{4n\pi}{m+2}-\delta.
\end{equation}
Thus, the second part of the theorem is proved by induction.

Next in order to  prove the first part of the theorem for the sector \eqref{zero1_sector}, we will determine which term in \eqref{asy_sector} dominates as $\lambda\to\infty$. To do that, we look at
\begin{align}
&\Re\left(L(G^{-1}(a),\omega^{-2}\lambda)-L(G^{n-1}(a),\omega^{2(n-1)}\lambda)\right)-\Re\left(L(G^{n}(a),\omega^{2n}\lambda)-L(a,\lambda)\right)\nonumber\\
&=K_m\left[\cos\left(-\frac{2\pi}{m}+\frac{m+2}{2m}\arg(\lambda)\right)-\cos\left(\frac{2(n-1)\pi}{m}+\frac{m+2}{2m}\arg(\lambda)\right)\right.\nonumber\\
&-\left.\left(\cos\left(-\frac{2n\pi}{m}+\frac{m+2}{2m}\arg(\lambda)\right)-\cos\left(\frac{m+2}{2m}\arg(\lambda)\right)\right)\right]|\lambda|^{\frac{1}{2}+\frac{1}{m}}(1+o(1))\nonumber\\
&=2K_m\sin\left(\frac{n\pi}{m}\right)\left[\sin\left(\frac{(n-2)\pi}{m}+\frac{m+2}{2m}\arg(\lambda)\right)\right.\nonumber\\
&\qquad\qquad\qquad\qquad\qquad\qquad\left.+\sin\left(\frac{n\pi}{m}+\frac{m+2}{2m}\arg(\lambda)\right)\right]|\lambda|^{\frac{1}{2}+\frac{1}{m}}(1+o(1))\nonumber\\
&=4K_m\sin\left(\frac{n\pi}{m}\right)\cos\left(\frac{\pi}{m}\right)\sin\left(\frac{(n-1)\pi}{m}+\frac{m+2}{2m}\arg(\lambda)\right)|\lambda|^{\frac{1}{2}+\frac{1}{m}}(1+o(1)),\label{dom_set}
\end{align}
that tends to positive infinity as $\lambda\to\infty$  (and hence the second term in \eqref{asy_sector} dominates) if $-\frac{2(n-1)\pi}{m+2}+\delta\leq\arg(\lambda)\leq \pi-\frac{4n\pi}{m+2}-\delta$.

We still need to prove \eqref{up_asy} for the sector 
\begin{equation}\label{asy_rem}
\pi-\frac{4n\pi}{m+2}-\delta\leq\arg(\lambda)\leq \pi-\frac{4n\pi}{m+2}+\delta,
\end{equation}
for which we use induction on $n$ again.

When $n=1$, \eqref{up_asy} holds by Lemma \ref{lemma_up}.

Suppose that \eqref{up_asy} in the sector \eqref{asy_rem} for $n-1$ with  $2\leq n<\lfloor\frac{m}{2}\rfloor$. 
Then \eqref{kplus1} and \eqref{kn_asy} with $k=n+1$  yield
\begin{equation}\nonumber
W_{-1,n}(a,\lambda)
=\frac{2\omega^{\mu(G^{-1}(a))}W_{0,1}(G^n(a),\omega^{2n}\lambda)}{\omega^{n-1}W_{0,n+1}(a,\lambda)}+\frac{W_{-1,n-1}(G(a),\omega^2\lambda)W_{0,n+2}(G^{-1}(a),\omega^{-2}\lambda)}{W_{0,n+1}(a,\lambda)}.
\end{equation}
If $\pi-\frac{4n\pi}{m+2}-\delta\leq\arg(\lambda)\leq \pi-\frac{4n\pi}{m+2}+\delta$, then $\pi-\frac{4(n-1)\pi}{m+2}-\delta\leq\arg(\omega^2\lambda)\leq \pi-\frac{4(n-1)\pi}{m+2}+\delta$. So
\begin{align}
W_{-1,n}(a,\lambda)
&=\frac{4\omega^{\mu(G^{-1}(a))+\mu(G^n(a))}}{\omega^{n-1}W_{0,n+1}(a,\lambda)}+\frac{W_{-1,n-1}(G(a),\omega^2\lambda)W_{0,n+2}(G^{-1}(a),\omega^{-2}\lambda)}{W_{0,n+1}(a,\lambda)}\label{last_eq}\\
&=\frac{4\omega^{\mu(G^{-1}(a))+\mu(G^n(a))}}{\omega^{n-1}[2i\omega^{-\frac{n+1}{2}}+o(1)]\exp\left[L(G^{n+1}(a),\omega^{2n-m}\lambda)+L(a,\lambda)\right]}\nonumber\\
&+\frac{[2\omega^{\frac{3-n}{2}+\mu(G^{n-1}(a))}+o(1)]\exp\left[L(a,\lambda)-L(G^{n-1}(a),\omega^{2(n-1)}\lambda)\right]}{[2i\omega^{-\frac{n+1}{2}}+o(1)]\exp\left[L(G^{n+1}(a),\omega^{2n-m}\lambda)+L(a,\lambda)\right]}\nonumber\\
&\times [2i\omega^{-\frac{n+2}{2}}+o(1)]\exp\left[L(G^{n+1}(a),\omega^{2n-m}\lambda)+L(G^{-1}(a),\omega^{-2}\lambda)\right]\nonumber\\
&=[-2i\omega^{\frac{3-n}{2}+\mu(G^{-1}(a))+\mu(G^n(a))}+o(1)]\exp\left[-L(G^{n+1}(a),\omega^{2n-m}\lambda)-L(a,\lambda)\right]\nonumber\\
&+[2\omega^{\frac{2-n}{2}+\mu(G^{n-1}(a))}+o(1)]\exp\left[L(G^{-1}(a),\omega^{-2}\lambda)-L(G^{n-1}(a),\omega^{2(n-1)}\lambda)\right],\nonumber
\end{align}
where we use the induction hypothesis for $W_{-1,n-1}(G(a),\omega^2\lambda)$ and use \eqref{sec_eq1} for $W_{0,n+1}(a,\lambda)$ and $W_{0,n+2}(G^{-1}(a),\omega^{-2}\lambda)$.
Next, we use an argument similar \eqref{dom_set} to complete the induction step. Thus, the theorem is proved.
\end{proof}
Next we investigate $W_{0,\lfloor\frac{m}{2}\rfloor}(a,\lambda)$.
\begin{theorem}\label{bod_thm}
If $m\geq 4$ is an even integer, then
\begin{align}
W_{-1,\lfloor\frac{m}{2}\rfloor}(a,\lambda)
&=-[2\omega^{2+\mu(G^{-1}(a))+\mu(G^{\frac{m}{2}}(a))}+o(1)]\exp\left[-L(G^{\frac{m+2}{2}}(a),\lambda)-L(a,\lambda)\right]\nonumber\\
&-[2\omega^{2+\mu(a)+\mu(G^{\frac{m-2}{2}}(a))}+o(1)]\exp\left[-L(G^{\frac{m-2}{2}}(a),\omega^{m-2}\lambda)-L(G^{m}(a),\omega^{m-2}\lambda)\right].\label{thm_eq1}
\end{align}
as $\lambda\to\infty$ in the sector 
\begin{equation}\label{zero4_sector}
-\pi+\frac{4\pi}{m+2}-\delta\leq\arg(\lambda)\leq-\pi+\frac{4\pi}{m+2}+\delta.
\end{equation}

If $m\geq 4$ is an odd integer, then
\begin{align}
W_{-1,\lfloor\frac{m}{2}\rfloor}(a,\lambda)
&=[2\omega^{\frac{5}{4}}+o(1)]\exp\left[-L(G^{\frac{m+1}{2}}(a),\omega^{-1}\lambda)-L(a,\lambda)\right]\nonumber\\
&+[2\omega^{\frac{5}{4}}+o(1)]\exp\left[L(G^{m+1}(a),\omega^{-2}\lambda)-L(G^{\frac{m-3}{2}}(a),\omega^{m-3}\lambda)\right].\label{thm_eq2}
\end{align}
as $\lambda\to\infty$ in the sector 
\begin{equation}\nonumber
-\pi+\frac{4\pi}{m+2}+\delta\leq\arg(\lambda)\leq-\pi+\frac{6\pi}{m+2}+\delta.
\end{equation}
\end{theorem}
\begin{proof}
We will use \eqref{last_eq} with $n=\lfloor\frac{m}{2}\rfloor$, that is,
\begin{equation}\nonumber
W_{-1,\lfloor\frac{m}{2}\rfloor}(a,\lambda)=\frac{4\omega^{\mu(G^{-1}(a))+\mu(G^{\lfloor\frac{m}{2}\rfloor}(a))}}{\omega^{\lfloor\frac{m}{2}\rfloor-1}W_{0,\lfloor\frac{m}{2}\rfloor+1}(a,\lambda)}+\frac{W_{-1,\lfloor\frac{m}{2}\rfloor-1}(G(a),\omega^2\lambda)W_{0,\lfloor\frac{m}{2}\rfloor+2}(G^{-1}(a),\omega^{-2}\lambda)}{W_{0,\lfloor\frac{m}{2}\rfloor+1}(a,\lambda)}.
\end{equation}
When $m$ is even, say $m=2k$, 
\begin{align}
W_{0,k+2}(G^{-1}(a),\omega^{-2}\lambda)&=W_{m+2,k+2}(G^{-1}(a),\omega^{-2}\lambda)\nonumber\\
&=-\omega^{-k-3}W_{-1,k-1}(G^{k+2}(a),\omega^{2k+4}\lambda)\nonumber\\
&=\omega^{-2}W_{-1,k-1}(G^{k+2}(a),\omega^{2}\lambda).\nonumber
\end{align}
So
\begin{align}
W_{-1,k}(a,\lambda)&=\frac{4\omega^{\mu(G^{-1}(a))+\mu(G^{k}(a))}}{\omega^{k-1}W_{0,k+1}(a,\lambda)}+\frac{W_{-1,k-1}(G(a),\omega^2\lambda)W_{0,k+2}(G^{-1}(a),\omega^{-2}\lambda)}{W_{0,k+1}(a,\lambda)}\nonumber\\
&=\frac{4\omega^{\mu(G^{-1}(a))+\mu(G^{k}(a))}}{\omega^{k-1}W_{0,k+1}(a,\lambda)}+\frac{W_{-1,k-1}(G(a),\omega^2\lambda)W_{-1,k-1}(G^{k+2}(a),\omega^{2}\lambda)}{\omega^2 W_{0,k+1}(a,\lambda)}\nonumber
\end{align}
Since $\lambda$ lies in \eqref{zero4_sector}, 
$$
-\frac{2\left(\lfloor\frac{m}{2}\rfloor-2\right)\pi}{m+2}+\delta\leq-\pi+\frac{8\pi}{m+2}-\delta\leq\arg(\omega^2\lambda)\leq\pi-\frac{4\left(\lfloor\frac{m}{2}\rfloor-1\right)\pi}{m+2}+\delta.
$$
\begin{align}
&W_{-1,k}(a,\lambda)\nonumber\\
&=\frac{4\omega^{\mu(G^{-1}(a))+\mu(G^{k}(a))}}{\omega^{k-1}W_{0,k+1}(a,\lambda)}+\frac{W_{-1,k-1}(G(a),\omega^2\lambda)W_{0,k+2}(G^{-1}(a),\omega^{-2}\lambda)}{W_{0,k+1}(a,\lambda)}\nonumber\\
&=\frac{4\omega^{\mu(G^{-1}(a))+\mu(G^{k}(a))}}{\omega^{k-1}W_{0,k+1}(a,\lambda)}+\frac{W_{-1,k-1}(G(a),\omega^2\lambda)W_{-1,k-1}(G^{k+2}(a),\omega^{2}\lambda)}{\omega^2 W_{0,k+1}(a,\lambda)}\nonumber\\
&=\frac{4\omega^{\mu(G^{-1}(a))+\mu(G^{k}(a))}}{\omega^{k-1}[2i\omega^{-\frac{k+1}{2}}+o(1)]\exp\left[L(G^{k+1}(a),\omega^{2k-m}\lambda)+L(a,\lambda)\right]}\nonumber\\
&+\frac{[2\omega^{\frac{3-k}{2}+\mu(G^{k-1}(a))}+o(1)]\exp\left[L(a,\lambda)-L(G^{k-1}(a),\omega^{2(k-1)}\lambda)\right]}{\omega^2 [2i\omega^{-\frac{k+1}{2}}+o(1)]\exp\left[L(G^{k+1}(a),\omega^{2k-m}\lambda)+L(a,\lambda)\right]}\nonumber\\
&\times[2\omega^{\frac{3-k}{2}+\mu(G^{2k}(a))}+o(1)]\exp\left[L(G^{k+1}(a),\lambda)-L(G^{2k}(a),\omega^{2(k-1)}\lambda)\right]\nonumber\\
&=[-2i\omega^{\frac{3-k}{2}+\mu(G^{-1}(a))+\mu(G^{k}(a))}+o(1)]\exp\left[-L(G^{k+1}(a),\lambda)-L(a,\lambda)\right]\nonumber\\
&-[2i\omega^{\frac{3-k}{2}+\mu(G^{k-1}(a))+\mu(G^{2k}(a))}+o(1)]\exp\left[-L(G^{k-1}(a),\omega^{2(k-1)}\lambda)-L(G^{2k}(a),\omega^{2(k-1)}\lambda)\right]\nonumber\\
&=-[2i\omega^{\frac{3-k}{2}+\mu(G^{-1}(a))+\mu(G^{k}(a))}+o(1)]\exp\left[-L(G^{k+1}(a),\lambda)-L(a,\lambda)\right]\nonumber\\
&-[2i\omega^{\frac{3-k}{2}+\mu(a)+\mu(G^{k-1}(a))}+o(1)]\exp\left[-L(G^{k-1}(a),\omega^{2(k-1)}\lambda)-L(G^{2k}(a),\omega^{2(k-1)}\lambda)\right],\nonumber
\end{align}
where we used \eqref{up_asy} for $W_{-1,k-1}(G(a),\cdot)$ and $W_{-1,k-1}(G^{k+2}(a),\cdot)$, and \eqref{sec_eq1} for $W_{0,k+1}(a,\cdot)$. Finally, we use $\omega^{-\frac{m+2}{4}}=-i$, to get the desired asymptotic expansion of $W_{-1,\lfloor\frac{m}{2}\rfloor}(a,\lambda)$.

Next we investigate the case when $m$ is odd, say $m=2k+1$ (so $\lfloor\frac{m}{2}\rfloor=k$).
\begin{equation}\nonumber
W_{0,\lfloor\frac{m}{2}\rfloor+1}(a,\lambda)=\omega^{-1}W_{-1,\lfloor\frac{m}{2}\rfloor}(G(a),\omega^{2}\lambda)
\end{equation}
and
\begin{align}
W_{0,k+2}(G^{-1}(a),\omega^{-2}\lambda)&=W_{m+2,k+2}(G^{-1}(a),\omega^{-2}\lambda)\nonumber\\
&=-\omega^{-k-2}W_{0,k+1}(G^{k+1}(a),\omega^{2k+2}\lambda)\nonumber\\
&=\omega^{-\frac{1}{2}}W_{0,k+1}(G^{k+1}(a),\omega^{-1}\lambda).\nonumber
\end{align}
Similarly to the proof of the theorem for $m$ even,
\begin{align}
W_{-1,k}(a,\lambda)
&=\frac{4\omega^{\mu(G^{-1}(a))+\mu(G^{k}(a))}}{\omega^{k-1}W_{0,k+1}(a,\lambda)}+\frac{W_{-1,k-1}(G(a),\omega^2\lambda)W_{0,k+2}(G^{-1}(a),\omega^{-2}\lambda)}{W_{0,k+1}(a,\lambda)}\nonumber\\
&=\frac{4\omega^{\mu(G^{-1}(a))+\mu(G^{k}(a))}}{\omega^{k-1}W_{0,k+1}(a,\lambda)}+\frac{W_{-1,k-1}(G(a),\omega^2\lambda)W_{0,k+1}(G^{k+1}(a),\omega^{-1}\lambda)}{\omega^{\frac{1}{2}} W_{0,k+1}(a,\lambda)}\nonumber\\
&=\frac{4\omega^{\mu(G^{-1}(a))+\mu(G^{k}(a))}}{\omega^{k-1}[2i\omega^{-\frac{k+1}{2}}+o(1)]\exp\left[L(G^{k+1}(a),\omega^{2k-m}\lambda)+L(a,\lambda)\right]}\nonumber\\
&+\frac{[2\omega^{\frac{3-k}{2}+\mu(G^{k-1}(a))}+o(1)]\exp\left[L(a,\lambda)-L(G^{k-1}(a),\omega^{2(k-1)}\lambda)\right]}{\omega^{\frac{1}{2}} [2i\omega^{-\frac{k+1}{2}}+o(1)]\exp\left[L(G^{k+1}(a),\omega^{2k-m}\lambda)+L(a,\lambda)\right]}\nonumber\\
&\times[2i\omega^{-\frac{k+1}{2}}+o(1)]\exp\left[L(G^{2k+2}(a),\omega^{2k-m-1}\lambda)+L(G^{k+1}(a),\omega^{-1}\lambda)\right]\nonumber\\
&=[-2i\omega^{\frac{3-k}{2}+\mu(G^{-1}(a))+\mu(G^{k}(a))}+o(1)]\exp\left[-L(G^{k+1}(a),\omega^{-1}\lambda)-L(a,\lambda)\right]\nonumber\\
&+[2\omega^{\frac{2-k}{2}+\mu(G^{k-1}(a))}+o(1)]\exp\left[L(G^{2k+2}(a),\omega^{-2}\lambda)-L(G^{k-1}(a),\omega^{2(k-1)}\lambda)\right]\nonumber\\
&=-[2i\omega^{\frac{3-k}{2}+\frac{m}{2}}+o(1)]\exp\left[-L(G^{k+1}(a),\omega^{-1}\lambda)-L(a,\lambda)\right]\nonumber\\
&+[2\omega^{\frac{2-k}{2}+\frac{m}{4}}+o(1)]\exp\left[L(G^{2k+2}(a),\omega^{-2}\lambda)-L(G^{k-1}(a),\omega^{2(k-1)}\lambda)\right],\nonumber
\end{align}
where we use \eqref{up_asy} for $W_{-1,k-1}(G(a),\cdot)$, and use \eqref{sec_eq1} for $W_{0,k+1}(a,\cdot)$ and  $W_{0,k+1}(G^{k+1}(a),\cdot)$. Finally, we use $\omega^{\frac{m+2}{4}}=i$, to get the asymptotic expansion of $W_{-1,\lfloor\frac{m}{2}\rfloor}(a,\lambda)$.
This  completes the proof.
\end{proof}

 The {\it order of an entire function} $g$ is defined by
$$\limsup_{r\rightarrow \infty}\frac{\log \log M(r,g)}{\log r},$$
where $M(r, g)=\max \{|g(re^{i\theta})|: 0\leq \theta\leq 2\pi\}$ for $r>0$.
If for some positive real numbers $\sigma,\, c_1,\, c_2$, we have $\exp[c_1 r^{\sigma}]\leq M(r,g)\leq  \exp[c_2 r^{\sigma}]$ for all large $r$, then the order of $g$ is $\sigma$. 
\begin{corollary}\label{zero_free}
Let $1\leq n\leq \frac{m}{2}$. Then
 the entire functions $W_{-1,n}(a,\cdot)$  are of order $\frac{1}{2}+\frac{1}{m}$, and hence they have infinitely many zeros in the complex plane. Moreover,  $W_{-1,n}(a,\cdot)$ have at most finitely many zeros in the sector 
$$
-\frac{2(n-1)\pi}{m+2}+\delta\leq\arg(\lambda)\leq\pi+\frac{4\pi}{m+2}-\delta.
$$
\end{corollary}

%%%%%%%%%%%%%%%%%%%%%%%%%%%%%%%%%%%%%%%%%%%%%%%%%%%
%%%%%%%%%%%%%%%%%%%%%%%%%%%%%%%%%%%%%%%%%%%%%%%%%%
\section{Relation between eigenvalues of $H_{\ell,P}$ and zeros of $W_{-1,n}(a,\cdot)$}\label{asymp_eigen}
%%%%%%%%%%%%%%%%%%%%%%%%%%%%%%%%%%%%%%%%%%%%%%%%%%%%
%%%%%%%%%%%%%%%%%%%%%%%%%%%%%%%%%%%%%%%%%%%%%%%%%%%%%%%
In this section, we will relate the eigenvalues of $H_{\ell,P}$ with zeros of some entire function $W_{-1,n}(a,\cdot)$.

Suppose $\ell=2k-1$ is odd with $1\leq \ell=2k-1\leq m-1$. Then \eqref{ptsym} becomes \eqref{rotated} by the scaling $v(z,\lambda)=u(-iz,\lambda)$, and $v(\cdot,\lambda)$ decays in the Stokes sectors $S_{-k}$ and $S_{k}$. Since $f_{k-1}$ and $f_{k}$ are linearly independent, for some $D_k$ and $\widetilde{D}_k$ one can write
$$f_{-k}(z,a,\lambda)=D_k(a,\lambda)f_{k-1}(z,a,\lambda)+\widetilde{D}_k(a,\lambda)f_{k}(z,a,\lambda).$$
Then one finds
$$
D_k(a,\lambda)=\frac{W_{-k,k}(a,\lambda)}{W_{k-1,k}(a,\lambda)}\quad\text{and}\quad \widetilde{D}_k(a,\lambda)=\frac{W_{-k,k-1}(a,\lambda)}{W_{k,k-1}(a,\lambda)}.
$$
Also it is easy to see that $\lambda$ is an eigenvalue of $H_{\ell,P}$ if and only if $D_k(a,\lambda)=0$ if and only if $W_{-k,k}(a,\lambda)=0$.
Since $W_{-k,k}(a,\lambda)=\omega^{k-1}W_{-1,2k-1}(G^{-k+1}(a),\omega^{-2k+2}\lambda)$, by Corollary \ref{zero_free} $W_{-k,k}(a,\lambda)$ has at most finitely many zeros in the sector 
$
-\frac{2(2k-2)\pi}{m+2}+\delta\leq\arg(\omega^{-2k+2}\lambda)\leq\pi+\frac{4\pi}{m+2}-\delta,
$
that is,
$$
\delta\leq\arg(\lambda)\leq\pi+\frac{4k\pi}{m+2}-\delta.
$$
This is true for all $a\in\C^{m-1}$ fixed. 

Next, by symmetry one can show that 
$W_{-k,k}(a,\lambda)$ has at most finitely many zeros in the sector $\pi\leq\arg(\lambda)\leq 2\pi-\delta$. We look at $H_{\ell,P}$ with $P(z)$ replaced by $\overline{P(\overline{z})}$ whose coefficient vector is $\overline{a}:=(\overline{a}_1,\,\overline{a}_2,\dots,\overline{a}_{m-1})$. Then one sees that
$W_{-k,k}(a,\lambda)=0$ if and only if $W_{-k,k}(\overline{a},\overline{\lambda})=0$. The latter has at most finitely many zeros in the sector $\delta\leq\arg(\overline{\lambda})\leq\pi.$  Thus, $W_{-k,k}(a,\lambda)$ has at most finitely many zeros in the sector $\pi\leq\arg(\lambda)\leq 2\pi-\delta$, and has infinitely many zeros in the sector $|\arg(\lambda)|\leq\delta$ since it is an entire function of order $\frac{1}{2}+\frac{1}{m}\in(0,1)$.

Suppose that $\ell=2k$ is even with $1\leq \ell=2k\leq m-1$. Then \eqref{ptsym} becomes \eqref{rotated2} by the scaling $y(z,\lambda)=u(-i\omega^{-\frac{1}{2}}z,\lambda)$, and $y(\cdot,\lambda)$ decays in the Stokes sectors $S_{-k}$ and $S_{k+1}$. We then see that the coefficient vector $\widetilde{a}$ of the polynomial $\omega^{-1}P(\omega^{-\frac{1}{2}}z)$ becomes
$$\widetilde{a}=G^{-\frac{1}{2}}(a).$$

Now one can express $f_{-k}$ as a linear combination of 
 $f_{k}$ and $f_{k+1}$ as follows. 
$$f_{-k}(z,\widetilde{a} ,\omega^{-1}\lambda)=\frac{W_{-k,k+1}(\widetilde{a},\omega^{-1}\lambda)}{W_{k,k+1}(\widetilde{a},\omega^{-1}\lambda)}f_{k}(z,\widetilde{a},\omega^{-1}\lambda)+\frac{W_{-k,k}(\widetilde{a},\omega^{-1}\lambda)}{W_{k+1,k}(\widetilde{a},\omega^{-1}\lambda)}f_{k+1}(z,\widetilde{a},\omega^{-1}\lambda).$$
Thus,  $\lambda$ is an eigenvalue of $H_{\ell,P}$ if and only if $W_{-k,k+1}(\widetilde{a},\omega^{-1}\lambda)=0$.
Since
$$W_{-k,k+1}(\widetilde{a},\omega^{-1}\lambda)=\omega^{k-1}W_{-1,2k}(G^{-k+1}(\widetilde{a}),\omega^{-2k+1}\lambda),$$ 
  by Corollary \ref{zero_free},
 $W_{-k,k+1}(\widetilde{a},\omega^{-1}\lambda)$ has at most finitely many zeros in the sector 
$
-\frac{2(2k-1)\pi}{m+2}+\delta\leq\arg(\omega^{-2k+1}\lambda)\leq\pi+\frac{4\pi}{m+2}-\delta,
$
that is,
$$
\delta\leq\arg(\lambda)\leq\pi+\frac{4k\pi}{m+2}-\delta.
$$
This is true for each $a\in\C^{m-1}$. So one can show that 
$W_{-k,k+1}(\widetilde{a},\omega^{-1}\lambda)$ has at most finitely many zeros in the sector $\pi\leq\arg(\lambda)\leq 2\pi-\delta$ by symmetry, similar to the cases when $\ell$ is odd.   Thus, $W_{-k,k+1}(a,\lambda)$  has infinitely many zeros in the sector $|\arg(\lambda)|\leq\delta$.

%%%%%%%%%%%%%%%%%%%%

\section{Proof of Theorem \ref{main_thm1} when $1\leq\ell<\lfloor\frac{m}{2}\rfloor$}\label{sec_7}
In this section, we prove Theorem \ref{main_thm1} for $1\leq\ell<\lfloor\frac{m}{2}\rfloor$. 

We first treat the case when $\ell$ is odd. 
\begin{proof}[Proof of Theorem ~\ref{main_thm1} when $1\leq\ell<\lfloor\frac{m}{2}\rfloor$ is odd]
Suppose that $1\leq \ell=2k-1<\lfloor\frac{m}{2}\rfloor$.
Recall that when $\ell$ is odd, $\lambda$ is an eigenvalue of $H_{\ell,P}$ if and only if 
 $W_{-k,k}(a,\lambda)=0$. Then since $W_{-k,k}(a,\lambda)$ has at most finitely many zeros outside the sector $|\arg(\lambda)|\leq\delta$, 
all the eigenvalues $\lambda$ of  $H_{\ell,P}$ lie in the sector $|\arg(\lambda)|\leq\delta$ if $|\lambda|$ is large enough.

Since 
\begin{equation}\nonumber
W_{-k,k}(a,\lambda)=\omega^{k-1}W_{-1,2k-1}(G^{-k+1}(a),\omega^{-2k+2}\lambda),
\end{equation}
 we will use \eqref{zero_asy} to investigate asymptotics of large eigenvalues. Suppose that $W_{-k,k}(a,\lambda)=0$ and $|\lambda|$ is large enough. Then from \eqref{zero_asy} with $n=2k-1$, and with $a$ and $\lambda$ replaced by $G^{-k+1}(a)$ and $\omega^{-2k+2}\lambda$, respectively, we have
\begin{align}
\left[1+o(1)\right]&\exp\left[L(G^{k}(a),\omega^{2k}\lambda)-L(G^{-k}(a),\omega^{-2k}\lambda)\right]\nonumber\\
&\times\exp\left[L(G^{k-1}(a),\omega^{2k-2}\lambda)-L(G^{-k+1}(a),\omega^{-2k+2}\lambda)\right]=-\omega^{2\nu(G^k(a))}.\nonumber
\end{align}
Next the term $\left[1+o(1)\right]$ can be  absorbed into the exponential function so that we get 
\begin{align}
&\exp\left[L(G^{k}(a),\omega^{2k}\lambda)-L(G^{-k}(a),\omega^{-2k}\lambda)\right]\nonumber\\
&\times\exp\left[L(G^{k-1}(a),\omega^{2k-2}\lambda)-L(G^{-k+1}(a),\omega^{-2k+2}\lambda)+o(1)\right]=-\omega^{2\nu(G^k(a))}.\label{sim_eq1}
\end{align}

For each odd integer $1\leq\ell=2k-1<\lfloor\frac{m}{2}\rfloor$, we define
\begin{align}
h_{m,\ell}(\lambda)=& L(G^{k}(a),\omega^{2k}\lambda)-L(G^{-k}(a),\omega^{-2k}\lambda)\nonumber\\
&+L(G^{k-1}(a),\omega^{2k-2}\lambda)-L(G^{-k+1}(a),\omega^{-2k+2}\lambda)+o(1),\nonumber
\end{align}
where the error term $o(1)$ is the same as that in \eqref{sim_eq1}.
Then by Corollary \ref{lemma_decay},
\begin{align}
h_{m,\ell}(\lambda)&=K_m\left(e^{\frac{2k\pi}{m}i}-e^{-\frac{2k\pi}{m}i}+e^{\frac{2(k-1)\pi}{m}i}-e^{-\frac{2(k-1)\pi}{m}i}\right)\lambda^{\frac{m+2}{2m}}(1+o(1))\nonumber\\
&=2iK_m\left(\sin\left(\frac{2k\pi}{m}\right)+\sin\left(\frac{2k\pi}{m}-\frac{2\pi}{m}\right)\right)\lambda^{\frac{m+2}{2m}}(1+o(1))\nonumber\\
&=4iK_m\cos\left(\frac{\pi}{m}\right)\sin\left(\frac{(2k-1)\pi}{m}\right)\lambda^{\frac{m+2}{2m}}(1+o(1))\quad\text{as}\quad\lambda\to\infty,\label{assco_eq}
\end{align}
in the sector $|\arg(\lambda)|\leq\delta$. Since $K_m>0$ and $0<\frac{(2k-1)\pi}{m}<\pi$, the function $h_{m,\ell}(\cdot)$ maps the region $|\lambda|\geq M_1$ for some large $M_1$ and $|\arg(\lambda)|\leq\delta$ into a region containing $|\lambda|\geq M_2$ for some large $M_2$ and $|\arg(\lambda)-\frac{\pi}{2}|\leq\varepsilon_1$ for some $\varepsilon_1>0$. So there exists a sequence of $\lambda_{n}$ in the sector $|\arg(\lambda)|\leq\delta$ such that $\exp\left[h_{m,\ell}(\lambda_{n})\right]= -\omega^{2\nu(G^k(a))}$ for all large $n\in\N$.

So, from \eqref{sim_eq1}, we have
\begin{equation}\nonumber
h_{m,\ell}(\lambda_{n})
=\ln\left(-\omega^{2\nu(G^k(a))}\right)
=\left(\frac{4\nu(G^k(a))}{m+2}+2n+1\right)\pi i,\quad\text{for all large $n\in\N$}.
\end{equation}
Thus, since 
\begin{align}
&-\frac{\nu(G^{k}(a))}{m}\ln(\omega^{2k}\lambda)+\frac{\nu(G^{-k}(a))}{m}\ln(\omega^{-2k}\lambda)\nonumber\\
&-\frac{\nu(G^{k-1}(a))}{m}\ln(\omega^{2k-2}\lambda)+\frac{\nu(G^{-k+1}(a))}{m}\ln(\omega^{-2k+2}\lambda)\nonumber\\
=& \frac{\nu(G^{k}(a))}{m}\left(-\ln(\omega^{2k}\lambda)+\ln(\omega^{-2k}\lambda)
+\ln(\omega^{2k-2}\lambda)-\ln(\omega^{-2k+2}\lambda)\right)\nonumber\\
=& -\frac{8\nu(G^{k}(a))}{m(m+2)}\pi i,\nonumber
\end{align}
from Lemma \ref{asy_lemma},
\begin{equation}\label{const_term}
\left(\frac{4\nu(G^k(a))}{m+2}+2n+1\right)\pi i=\sum_{j=0}^{\lfloor\frac{m+2}{2}\rfloor}d_{\ell,j}(a)\lambda_{n}^{\frac{1}{2}+\frac{1-j}{m}}-\frac{8\nu(G^k(a))}{m(m+2)}\pi i+o(1),
\end{equation}
where for $0\leq j\leq \frac{m+2}{2}$, the coefficients $d_{\ell,j}(a)$ are given by
\begin{align}
d_{\ell,j}(a)&=K_{m,j}(G^{k}(a))\omega^{2k\left(\frac{1}{2}+\frac{1-j}{m}\right)}-K_{m,j}(G^{-k}(a))\omega^{-2k\left(\frac{1}{2}+\frac{1-j}{m}\right)}\nonumber\\
&+K_{m,j}(G^{k-1}(a))\omega^{2(k-1)\left(\frac{1}{2}+\frac{1-j}{m}\right)}-K_{m,j}(G^{-k+1}(a))\omega^{-2(k-1)\left(\frac{1}{2}+\frac{1-j}{m}\right)}.\nonumber
\end{align}
Since $\ell=2k-1$, for $1\leq j\leq\frac{m+2}{2}$,
\begin{align}
d_{\ell,j}(a)&=K_{m,j}(G^{\frac{\ell+1}{2}}(a))\omega^{(\ell+1)\left(\frac{1}{2}+\frac{1-j}{m}\right)}-K_{m,j}(G^{-\frac{\ell+1}{2}}(a))\omega^{-(\ell+1)\left(\frac{1}{2}+\frac{1-j}{m}\right)}\nonumber\\
&+K_{m,j}(G^{\frac{\ell-1}{2}}(a))\omega^{(\ell-1)\left(\frac{1}{2}+\frac{1-j}{m}\right)}-K_{m,j}(G^{-\frac{\ell-1}{2}}(a))\omega^{-(\ell-1)\left(\frac{1}{2}+\frac{1-j}{m}\right)}\nonumber\\
&=4i\sum_{k=1}^j(-1)^kK_{m,j,k}b_{j,k}(a)\sin\left(\frac{(1-j)\ell\pi}{m}\right)\cos\left(\frac{(1-j)\pi}{m}\right),\nonumber
\end{align}
where we used Lemma \ref{lemma_25} along with \eqref{K_deF}.
Notice that \eqref{assco_eq} shows 
$$d_{\ell,0}(a)=4iK_{m,0}\sin\left(\frac{\ell\pi}{m}\right)\cos\left(\frac{\pi}{m}\right).$$
If $m$ is even, then we redefine $d_{\ell,\frac{m+2}{2}}(a)$ as the sum of $d_{\ell,\frac{m+2}{2}}(a)=0$ above and $-\frac{4\nu(G^k(a))}{m}\pi i$ (c. f., equation \eqref{const_term}).
Finally, note that $\nu(G^k(a))=(-1)^k\nu(a)=\eta_{m,\ell}(a).$ This completes the proof.
\end{proof}

Next we prove Theorem ~\ref{main_thm1} for $1\leq\ell<\lfloor\frac{m}{2}\rfloor$ is even.
\begin{proof}[Proof of Theorem ~\ref{main_thm1} when $1\leq\ell<\lfloor\frac{m}{2}\rfloor$ is even]

Proof is very similar to the case when $1\leq\ell<\lfloor\frac{m}{2}\rfloor$ is even. Let $\ell=2k$ for some $k\in\N$. 

Recall that  $\lambda$ is an eigenvalue of $H_{\ell,P}$ if and only if  $W_{-1,2k}(G^{-k+1}(\widetilde{a}),\omega^{-2k+1}\lambda)=0$. Then from \eqref{zero_asy},  we have
\begin{align}
&\exp\left[L(G^{k+1}(\widetilde{a}),\omega^{2k+1}\lambda)-L(G^{-k}(\widetilde{a}),\omega^{-2k-1}\lambda)\right]\nonumber\\
&\times\exp\left[L(G^{k}(\widetilde{a}),\omega^{2k-1}\lambda)-L(G^{-k+1}(\widetilde{a}),\omega^{-2k+1}\lambda)+o(1)\right]=-1,\label{sim_eq2}
\end{align}
where we used $\left[1+o(1)\right]=\exp[o(1)]$ again.

Like in the case when $1\leq\ell<\lfloor\frac{m}{2}\rfloor$ is odd, from Lemma \ref{asy_lemma},
\begin{equation}\nonumber
\left(2n+1\right)\pi i=\sum_{j=0}^{\frac{m+2}{2}}d_{\ell,j}(a)\lambda_{n}^{\frac{1}{2}+\frac{1-j}{m}}+o(1),
\end{equation}
where for $0\leq j\leq \frac{m+2}{2}$, 
\begin{align}
d_{\ell,j}(a)&=K_{m,j}(G^{k+1}(\widetilde{a}))\omega^{(2k+1)\left(\frac{1}{2}+\frac{1-j}{m}\right)}-K_{m,j}(G^{-k}(\widetilde{a}))\omega^{-(2k+1)\left(\frac{1}{2}+\frac{1-j}{m}\right)}\nonumber\\
&+K_{m,j}(G^{k}(\widetilde{a}))\omega^{(2k-1)\left(\frac{1}{2}+\frac{1-j}{m}\right)}-K_{m,j}(G^{-k+1}(\widetilde{a}))\omega^{-(2k-1)\left(\frac{1}{2}+\frac{1-j}{m}\right)}\nonumber\\
&=\sum_{k=1}^j(-1)^kK_{m,j,k}b_{j,k}(a)\left(\omega^{-j\frac{\ell+1}{2}+(\ell+1)\left(\frac{1}{2}+\frac{1-j}{m}\right)}-\omega^{j\frac{\ell+1}{2}-(\ell+1)\left(\frac{1}{2}+\frac{1-j}{m}\right)}\right.\nonumber\\
&\qquad\qquad\qquad\qquad\qquad\qquad\left.+\omega^{-j\frac{\ell-1}{2}+(\ell-1)\left(\frac{1}{2}+\frac{1-j}{m}\right)}-\omega^{j\frac{\ell-1}{2}-(\ell-1)\left(\frac{1}{2}+\frac{1-j}{m}\right)}\right)\nonumber\\
&=4i\sum_{k=1}^j(-1)^kK_{m,j,k}b_{j,k}(a)\sin\left(\frac{(1-j)\ell\pi}{m}\right)\cos\left(\frac{(1-j)\pi}{m}\right),\nonumber
\end{align}
where we used $\widetilde{a}=G^{-\frac{1}{2}}(a)$ and $\ell=2k$ as well.
\end{proof}

%%%%%%%%%%%%%%%%%%%%%
%%%%%%%%%%%%%%%%%%%%%%
\section{Proof of Theorem \ref{main_thm1} when $\ell=\lfloor\frac{m}{2}\rfloor$ and when $\frac{m}{2}<\ell\leq m-1$.}\label{sec_8}
In this section, we prove Theorem \ref{main_thm1} for $\ell=\lfloor\frac{m}{2}\rfloor$. We first prove the theorem when $m$ is even, and later, we will treat the cases when $m$ is odd. Then at the end of the section, we will prove the theorem when $\frac{m}{2}<\ell\leq m-1$, by scaling.
\subsection{When $m$ is even} We further divide the case into when $\ell$ is odd and when $\ell$ is even.
\begin{proof}[Proof of Theorem ~\ref{main_thm1} when $m$ is even and $\ell=\lfloor\frac{m}{2}\rfloor$ is odd]
Suppose that $m$ is even and $\ell=\frac{m}{2}=2k-1$ for some $k\in\N$. Since $\lambda$ is an eigenvalue of $H_{\ell,P}$ if and only if $W_{-k,k}(a,\lambda)=0$ if and only if  $W_{-1,2k-1}(G^{-k+1}(a),\omega^{-2k+2}\lambda)=0$.

Suppose that $W_{-1,2k-1}(G^{-k+1}(a),\omega^{-2k+2}\lambda)=0$. If $|\lambda|$ is large enough, then from  \eqref{thm_eq1},
\begin{align}
&\exp\left[L(G^{3k-1}(a),\omega^{2k-2}\lambda)-L(G^{k+1}(a),\omega^{-2k+2}\lambda)\right]\nonumber\\
&\times\exp\left[L(G^{k-1}(a),\omega^{2k-2}\lambda)-L(G^{-k+1}(a),\omega^{-2k+2}\lambda)+o(1)\right]=-\omega^{4\nu(G^{k}(a))}.\nonumber
\end{align}

From Lemma \ref{asy_lemma},
\begin{equation}\nonumber
\left(\frac{8\nu(G^{k}(a))}{m+2}+2n+1\right)\pi i=\sum_{j=0}^{\frac{m+2}{2}}d_{\ell,j}(a)\lambda_{n}^{\frac{1}{2}+\frac{1-j}{m}}+\frac{16(k-1)\nu(G^{k}(a))}{m(m+2)}\pi i+o(1),
\end{equation}
where for $0\leq j\leq \frac{m+2}{2}$, the coefficients $d_{\ell,j}(a)$ are given by
\begin{align}
d_{\ell,j}(a)&=K_{m,j}(G^{3k-1}(a))\omega^{2(k-1)\left(\frac{1}{2}+\frac{1-j}{m}\right)}-K_{m,j}(G^{k+1}(a))\omega^{-2(k-1)\left(\frac{1}{2}+\frac{1-j}{m}\right)}\nonumber\\
&+K_{m,j}(G^{k-1}(a))\omega^{2(k-1)\left(\frac{1}{2}+\frac{1-j}{m}\right)}-K_{m,j}(G^{-k+1}(a))\omega^{-2(k-1)\left(\frac{1}{2}+\frac{1-j}{m}\right)}\nonumber\\
&=K_{m,j}(G^{\frac{3m+2}{4}}(a))\omega^{\frac{m-2}{2}\left(\frac{1}{2}+\frac{1-j}{m}\right)}-K_{m,j}(G^{\frac{m+6}{4}}(a))\omega^{-\frac{m-2}{2}\left(\frac{1}{2}+\frac{1-j}{m}\right)}\nonumber\\
&+K_{m,j}(G^{\frac{m-2}{4}}(a))\omega^{\frac{m-2}{2}\left(\frac{1}{2}+\frac{1-j}{m}\right)}-K_{m,j}(G^{-\frac{m-2}{4}}(a))\omega^{-\frac{m-2}{2}\left(\frac{1}{2}+\frac{1-j}{m}\right)}\nonumber\\
&=K_{m,j}(a)\left(\omega^{-j\frac{3m+2}{4}+\frac{m-2}{2}\left(\frac{1}{2}+\frac{1-j}{m}\right)}-\omega^{-j\frac{m+6}{4}-\frac{m-2}{2}\left(\frac{1}{2}+\frac{1-j}{m}\right)}\right.\nonumber\\
&\qquad\qquad\qquad\qquad\qquad\qquad\left.\omega^{-j\frac{m-2}{4}+\frac{m-2}{2}\left(\frac{1}{2}+\frac{1-j}{m}\right)}-\omega^{j\frac{m-2}{4}-\frac{m-2}{2}\left(\frac{1}{2}+\frac{1-j}{m}\right)}\right)\nonumber\\
&=4iK_{m,j}(a)\sin\left(\frac{(1-j)\pi}{2}\right)\cos\left(\frac{(1-j)\pi}{m}\right),\nonumber
\end{align}
where we used Lemma \ref{lemma_25} as well as some other thing as before. In addition, we used 
\begin{align}
&\omega^{-j\frac{3m+2}{4}+\frac{m-2}{2}\left(\frac{1}{2}+\frac{1-j}{m}\right)}=-\omega^{j\frac{m+2}{4}-\frac{m+2}{2}\left(\frac{1}{2}+\frac{1-j}{m}\right)},\nonumber\\
&\omega^{-j\frac{m+6}{4}-\frac{m-2}{2}\left(\frac{1}{2}+\frac{1-j}{m}\right)}=-\omega^{-j\frac{m+2}{4}+\frac{m+2}{2}\left(\frac{1}{2}+\frac{1-j}{m}\right)}.\label{sine_eq}
\end{align}
\end{proof}

\begin{proof}[Proof of Theorem ~\ref{main_thm1} when $m$ is even and $\ell=\lfloor\frac{m}{2}\rfloor$ is even]
Let $\ell=\frac{m}{2}=2k$ for some $k\in\N$. Then
$\lambda$ is an eigenvalue of $H_{\ell,P}$ if and only if $W_{-1,2k}(G^{-k+1}(\widetilde{a}),\omega^{-2k+1}\lambda)=0$.
If $W_{-1,2k}(G^{-k+1}(\widetilde{a}),\omega^{-2k+1}\lambda)=0$ and $|\lambda|$ large enough, then from  \eqref{thm_eq1},
\begin{align}
&\exp\left[L(G^{3k+1}(\widetilde{a}),\omega^{2k-1}\lambda)-L((G^{k+2}(\widetilde{a}),\omega^{-2k+1}\lambda)\right]\nonumber\\
&\times\exp\left[L(G^{k}(\widetilde{a}),\omega^{2k-1}\lambda)-L(G^{-k+1}(\widetilde{a}),\omega^{-2k+1}\lambda)+o(1)\right]=-1.\nonumber
\end{align}

Then like before,
\begin{equation}\nonumber
\left(2n+1\right)\pi i=\sum_{j=0}^{\frac{m+2}{2}}d_{\ell,j}(\widetilde{a})\lambda_{n}^{\frac{1}{2}+\frac{1-j}{m}}+o(1),
\end{equation}
where for $0\leq j\leq \frac{m+2}{2}$, 
\begin{align}
d_{\ell,j}(a)&=K_{m,j}(G^{3k+1}(\widetilde{a}))\omega^{(2k-1)\left(\frac{1}{2}+\frac{1-j}{m}\right)}-K_{m,j}((G^{k+2}(\widetilde{a}))\omega^{-(2k-1)\left(\frac{1}{2}+\frac{1-j}{m}\right)}\nonumber\\
&+K_{m,j}(G^{k}(\widetilde{a}))\omega^{(2k-1)\left(\frac{1}{2}+\frac{1-j}{m}\right)}-K_{m,j}(G^{-k+1}(\widetilde{a}))\omega^{-(2k-1)\left(\frac{1}{2}+\frac{1-j}{m}\right)}\nonumber\\
&=K_{m,j}(G^{\frac{3m+2}{4}}(a))\omega^{\frac{m-2}{2}\left(\frac{1}{2}+\frac{1-j}{m}\right)}-K_{m,j}((G^{\frac{m+6}{4}}(a))\omega^{-\frac{m-2}{2}\left(\frac{1}{2}+\frac{1-j}{m}\right)}\nonumber\\
&+K_{m,j}(G^{\frac{m-2}{4}}(a))\omega^{\frac{m-2}{2}\left(\frac{1}{2}+\frac{1-j}{m}\right)}-K_{m,j}(G^{-\frac{m-2}{4}}(a))\omega^{-\frac{m-2}{2}\left(\frac{1}{2}+\frac{1-j}{m}\right)}\nonumber\\
&=4i\sum_{k=1}^j(-1)^kK_{m,j,k}b_{j,k}(a)\sin\left(\frac{(1-j)\pi}{2}\right)\cos\left(\frac{(1-j)\pi}{m}\right).\nonumber
\end{align}
\end{proof}

\subsection{When $m$ is odd} We divide the case into when $\ell$ is odd and when $\ell$ is even.
\begin{proof}[Proof of Theorem ~\ref{main_thm1} when $m$ and $\ell=\lfloor\frac{m}{2}\rfloor$ are odd]
Let $m$ and $\ell=\frac{m-1}{2}=2k-1$ be odd. Suppose that $\lambda$ is an eigenvalue of $H_{\ell,P}$. Since $\lambda$ is an eigenvalue of $H_{\ell,P}$ if and only if $W_{-1,2k-1}(G^{-k+1}(a),\omega^{-2k+2}\lambda)=0$,
 if  $|\lambda|$ is large enough, then from  \eqref{thm_eq2},
\begin{align}
&\exp\left[L(G^{k-1}(a),\omega^{2k-2}\lambda)-L(G^{3k+1}(a),\omega^{-2k}\lambda)\right]\nonumber\\
&\times\exp\left[-L(G^{k+1}(a),\omega^{-2k+1}\lambda)-L(G^{-k+1}(a),\omega^{-2k+2}\lambda)+o(1)\right]=-1.\nonumber
\end{align}

Then, 
from Lemma \ref{asy_lemma},
\begin{equation}\nonumber
\left(2n+1\right)\pi i=\sum_{j=0}^{\frac{m+1}{2}}d_{\ell,j}(a)\lambda_{n}^{\frac{1}{2}+\frac{1-j}{m}}+o(1),
\end{equation}
where for $0\leq j\leq \frac{m+2}{2}$, 
\begin{align}
d_{\ell,j}(a)&=K_{m,j}(G^{\frac{m-3}{4}}(a))\omega^{\frac{m-3}{2}\left(\frac{1}{2}+\frac{1-j}{m}\right)}-K_{m,j}(G^{\frac{3m+7}{4}}(a))\omega^{-\frac{m+1}{2}\left(\frac{1}{2}+\frac{1-j}{m}\right)}\nonumber\\
&-K_{m,j}(G^{\frac{m+5}{4}}(a))\omega^{-\frac{m-1}{2}\left(\frac{1}{2}+\frac{1-j}{m}\right)}-K_{m,j}(G^{-\frac{m-3}{4}}(a))\omega^{-\frac{m-3}{2}\left(\frac{1}{2}+\frac{1-j}{m}\right)}\nonumber\\
&=K_{m,j}(a)\left(\omega^{-j\frac{m-3}{4}+\frac{m-3}{2}\left(\frac{1}{2}+\frac{1-j}{m}\right)}-\omega^{-j\frac{3m+7}{4}-\frac{m+1}{2}\left(\frac{1}{2}+\frac{1-j}{m}\right)}\right.\nonumber\\
&\left.-\omega^{-j\frac{m+5}{4}-\frac{m-1}{2}\left(\frac{1}{2}+\frac{1-j}{m}\right)}-\omega^{j\frac{m-3}{4}-\frac{m-3}{2}\left(\frac{1}{2}+\frac{1-j}{m}\right)}\right)\nonumber\\
&=2iK_{m,j}(a)\left(\sin\left(\frac{(m-3)(1-j)\pi}{2m}\right)+\sin\left(\frac{(m+1)(1-j)\pi}{2m}\right)\right)\nonumber\\
&=2iK_{m,j}(a)\sin\left(\frac{(m-1)(1-j)\pi}{2m}\right)\cos\left(\frac{(1-j)\pi}{m}\right),\nonumber
\end{align}
where we used 
\begin{align}
&\omega^{-j\frac{3m+7}{4}-\frac{m+1}{2}\left(\frac{1}{2}+\frac{1-j}{m}\right)}=\omega^{j\frac{m+1}{4}-\frac{m+1}{2}\left(\frac{1}{2}+\frac{1-j}{m}\right)},\nonumber\\
& \omega^{-j\frac{m+5}{4}-\frac{m-1}{2}\left(\frac{1}{2}+\frac{1-j}{m}\right)}=-\omega^{-j\frac{m+1}{4}+\frac{m+1}{2}\left(\frac{1}{2}+\frac{1-j}{m}\right)}.\label{sine_eq2}
\end{align}
\end{proof}

\begin{proof}[Proof of Theorem ~\ref{main_thm1} when $m$ is odd and $\ell=\lfloor\frac{m}{2}\rfloor$ is even]
Let $\ell=\frac{m-1}{2}=2k$ for some $k\in\N$. Then
$\lambda$ is an eigenvalue of $H_{\ell,P}$ if and only if $W_{-1,2k}(G^{-k+1}(\widetilde{a}),\omega^{-2k+1}\lambda)=0$.
If $W_{-1,2k}(G^{-k+1}(\widetilde{a}),\omega^{-2k+1}\lambda)=0$ and $|\lambda|$ large enough, then from  \eqref{thm_eq2},
\begin{align}
&\exp\left[L(G^{k}(\widetilde{a}),\omega^{2k-1}\lambda)-L(G^{-k+1}(\widetilde{a}),\omega^{-2k+1}\lambda)\right]\nonumber\\
&\times\exp\left[-L(G^{3k+3}(\widetilde{a}),\omega^{-2k-1}\lambda)-L(G^{k+2}(\widetilde{a}),\omega^{-2k}\lambda)+o(1)\right]=-1.\nonumber
\end{align}

Then, 
from Lemma \ref{asy_lemma},
\begin{equation}\nonumber
\left(2n+1\right)\pi i=\sum_{j=0}^{\frac{m+1}{2}}d_{\ell,j}(\widetilde{a})\lambda_{n}^{\frac{1}{2}+\frac{1-j}{m}}+o(1),
\end{equation}
where for $0\leq j\leq \frac{m+2}{2}$, 
\begin{align}
d_{\ell,j}(a)&=K_{m,j}(G^{k}(\widetilde{a}))\omega^{(2k-1)\left(\frac{1}{2}+\frac{1-j}{m}\right)}-K_{m,j}((G^{-k+1}(\widetilde{a}))\omega^{-(2k-1)\left(\frac{1}{2}+\frac{1-j}{m}\right)}\nonumber\\
&-K_{m,j}(G^{3k+3}(\widetilde{a}))\omega^{-(2k+1)\left(\frac{1}{2}+\frac{1-j}{m}\right)}-K_{m,j}(G^{k+2}(\widetilde{a}))\omega^{-2k\left(\frac{1}{2}+\frac{1-j}{m}\right)}\nonumber\\
&=K_{m,j}(G^{\frac{m-3}{4}}(a))\omega^{\frac{m-3}{2}\left(\frac{1}{2}+\frac{1-j}{m}\right)}-K_{m,j}((G^{-\frac{m-3}{4}}(a))\omega^{-\frac{m-3}{2}\left(\frac{1}{2}+\frac{1-j}{m}\right)}\nonumber\\
&-K_{m,j}(G^{\frac{3m+7}{4}}(a))\omega^{-\frac{m+1}{2}\left(\frac{1}{2}+\frac{1-j}{m}\right)}-K_{m,j}(G^{\frac{m+5}{4}}(a))\omega^{-\frac{m-1}{2}\left(\frac{1}{2}+\frac{1-j}{m}\right)}\nonumber\\
&=\sum_{k=1}^j(-1)^kK_{m,j,k}b_{j,k}(a)\sin\left(\frac{(m-1)(1-j)\pi}{2m}\right)\cos\left(\frac{(1-j)\pi}{m}\right).\nonumber
\end{align}
\end{proof}
Theorem \ref{main_thm1} for $1\leq\ell\leq\frac{m}{2}$ has been proved. Next we prove Theorem \ref{main_thm1} for $\frac{m}{2}<\ell\leq m-1$, by scaling $z\mapsto -z$.
\subsection{When $\frac{m}{2}<\ell\leq m-1$}
\begin{proof}[Proof of Theorem ~\ref{main_thm1} when $\frac{m}{2}<\ell\leq m-1$]
Suppose that  $\ell>\frac{m}{2}$. If $u(\cdot,\lambda)$ is an eigenfunction of $H_{\ell,P}$, then $v(z,\lambda)=u(-z,\lambda)$ solves 
\begin{equation}\nonumber
-v^\dd(z,\lambda)+\left[(-1)^{-\ell}(-iz)^m-P(-iz)\right]v(z,\lambda)
=\lambda v(z,\lambda),
\end{equation}
and
\begin{equation}\nonumber
\text{$v(z)\rightarrow 0$ exponentially, as $z\rightarrow \infty$ along the two rays}\quad \arg z=-\frac{\pi}{2}\pm \frac{((m-\ell)+1)\pi}{m+2}.
\end{equation}
The coefficient vector of $P(-z)$ is
$((-1)^{m-1}a_1,(-1)^{m-2}a_2,\dots,(-1)^{1}a_{m-1}).$ Certainly, 
$$\sin\left(\frac{(m-\ell)\pi}{m}\right)=\sin\left(\frac{\ell\pi}{m}\right),\quad \sin\left(\frac{(1-j)(m-\ell)\pi}{m}\right)=(-1)^j\sin\left(\frac{(1-j)\ell\pi}{m}\right).$$ Also, one can find from \eqref{bjk_def} that for $1\leq k\leq j\leq\frac{m+2}{2}$,
$$b_{j,k}((-1)^{m-1}a_1,(-1)^{m-2}a_2,\dots,-a_{m-1})=(-1)^{mk-j}b_{j,k}(a_1,a_2,\dots,a_{m-1}).$$
Moreover,
$d_{m-\ell,j}((-1)^{m-1}a_1,(-1)^{m-2}a_2,\dots,-a_{m-1})=d_{\ell,j}(a_1,a_2,\dots,a_{m-1}).$
This completes proof of Theorem \ref{main_thm1}.
\end{proof}
%%%%%%%%%%%%%%%%%%%%%%%%%%%%%%%%%%%%%%%%%%%%%%%%%%%%%
%%%%%%%%%%%%%%%%%%%%%%%%%%%%%%%%%%%%%%%%%%%%%%%%%%%%%
\section{Proof of Theorem \ref{eigen_asy}}\label{sec_9}
%%%%%%%%%%%%%%%%%%%%%%%%%%%%%%%%%%%%%%%%%%%%%%%%%%%%%
%%%%%%%%%%%%%%%%%%%%%%%%%%%%%%%%%%%%%%%%%%%%%%%%%%%%
In this section, we prove Theorem \ref{eigen_asy}.

We begin with
\begin{equation}\nonumber
\left(2n+1\right)\pi i=\sum_{j=0}^{\lfloor\frac{m+2}{2}\rfloor}d_{\ell,j}(a)\lambda_{n}^{\frac{1}{2}+\frac{1-j}{m}}+o(1),\quad\text{as $n\to\infty$},
\end{equation}
\begin{equation}\nonumber
\frac{\left(2n+1\right)\pi i}{d_{\ell,0}(a)}=\lambda_{n}^{\frac{1}{2}+\frac{1}{m}}+\sum_{j=1}^{\lfloor\frac{m}{2}+1\rfloor}\frac{d_{\ell,j}(a)}{d_{\ell,0}(a)}\lambda_{n}^{\frac{1}{2}+\frac{1-j}{m}}+o(1).
\end{equation}

Let 
$$c_{j}(a)=\frac{d_{\ell,j}(a)}{d_{\ell,0}(a)},\quad 1\leq j\leq\frac{m+2}{2}.$$
Introduce the decomposition
$
\lambda_{n}=\lambda_{n,0}+\lambda_{n,1},
$
where 
\begin{equation}\nonumber
 \lambda_{n,0}=\left(\frac{\left(2n+1\right)\pi i}{d_{\ell,0}(a)}\right)^{\frac{2m}{m+2}}\,\,\text{and}\quad \frac{\lambda_{n,0}}{\lambda_{n,0}}=o\left(1\right).
\end{equation}
Then  we have
\begin{align}
\lambda_{n,0}^{\frac{1}{2}+\frac{1}{m}}
&=\lambda_{n,0}^{\frac{1}{2}+\frac{1}{m}}\left(1+\frac{\lambda_{n,1}}{\lambda_{n,0}}\right)^{\frac{1}{2}+\frac{1}{m}}+\sum_{j=1}^{\lfloor\frac{m}{2}+1\rfloor}c_{j}(a)\lambda_{n,0}^{\frac{1}{2}+\frac{1-j}{m}}\left(1+\frac{\lambda_{n,1}}{\lambda_{n,0}}\right)^{\frac{1}{2}+\frac{1-j}{m}}+o(1)\nonumber\\
&=\lambda_{n,0}^{\frac{1}{2}+\frac{1}{m}}\left(1+\sum_{k=1}^{\infty}{\frac{1}{2}+\frac{1}{m}\choose k}\left(\frac{\lambda_{n,1}}{\lambda_{n,0}}\right)^k\right)\nonumber\\
&+\sum_{j=1}^{\lfloor\frac{m}{2}+1\rfloor}c_{j}(a)\lambda_{n,0}^{\frac{1}{2}+\frac{1-j}{m}}\left(1+\sum_{k=1}^{\infty}{\frac{1}{2}+\frac{1-j}{m}\choose k}\left(\frac{\lambda_{n,1}}{\lambda_{n,0}}\right)^k\right)+o(1).\nonumber
\end{align}
Thus,
\begin{align}
0&={\frac{1}{2}+\frac{1}{m}\choose 1}\frac{\lambda_{n,1}}{\lambda_{n,0}}+\sum_{k=2}^{\infty}{\frac{1}{2}+\frac{1}{m}\choose k}\left(\frac{\lambda_{n,1}}{\lambda_{n,0}}\right)^k\nonumber\\
&+\sum_{j=1}^{\lfloor\frac{m}{2}+1\rfloor}c_{j}(a)\lambda_{n,0}^{-\frac{j}{m}}\left(1+\sum_{k=1}^{\infty}{\frac{1}{2}+\frac{1-j}{m}\choose k}\left(\frac{\lambda_{n,1}}{\lambda_{n,0}}\right)^k\right)+o\left(\lambda_{n,0}^{-\frac{1}{2}-\frac{1}{m}}\right),\nonumber
\end{align}
and hence
\begin{align}
&\frac{\lambda_{n,1}}{\lambda_{n,0}}+\sum_{k=2}^{\infty}\frac{{\frac{1}{2}+\frac{1}{m}\choose k}}{{\frac{1}{2}+\frac{1}{m}\choose 1}}\left(\frac{\lambda_{n,1}}{\lambda_{n,0}}\right)^k\nonumber\\
&+\sum_{j=1}^{\lfloor\frac{m}{2}+1\rfloor}c_{j}(a)\lambda_{n,0}^{-\frac{j}{m}}\left(\sum_{k=1}^{\infty}\frac{{\frac{1}{2}+\frac{1-j}{m}\choose k}}{{\frac{1}{2}+\frac{1}{m}\choose 1}}\left(\frac{\lambda_{n,1}}{\lambda_{n,0}}\right)^k\right)+o\left(\lambda_{n,0}^{-\frac{1}{2}-\frac{1}{m}}\right)\nonumber\\
&=-\frac{1}{{\frac{1}{2}+\frac{1}{m}\choose 1}}\sum_{j=1}^{\lfloor\frac{m}{2}+1\rfloor}c_{j}(a)\lambda_{n,0}^{-\frac{j}{m}}.\label{asy_eq4}
\end{align}
Thus, one concludes
$
\frac{\lambda_{n,1}}{\lambda_{n,0}}=\lambda_{n,2}+\lambda_{n,3},
$
where
\begin{equation}\label{ex_eq1}
\lambda_{n,2}=-\frac{1}{{\frac{1}{2}+\frac{1}{m}\choose 1}}c_{1}(a)\lambda_{n,0}^{-\frac{1}{m}}\,\,\text{ and }\,\,\lambda_{n,3}=o\left(\lambda_{n,0}^{-\frac{1}{m}}\right).
\end{equation}

Notice that $\lambda_{n,2}=0$ since $c_{1}(a)=0$. Hence, from  \eqref{ex_eq1} along with \eqref{asy_eq4} we have
\begin{align}
&\lambda_{n,3}+\sum_{{k_1}=2}^{\infty}\frac{{\frac{1}{2}+\frac{1}{m}\choose {k_1}}}{{\frac{1}{2}+\frac{1}{m}\choose 1}}\lambda_{n,3}^{k_1}+\sum_{j=2}^{\lfloor\frac{m}{2}+1\rfloor}c_{j}(a)\lambda_{n,0}^{-\frac{j}{m}}\left(\sum_{{k_1}=1}^{\infty}\frac{{\frac{1}{2}+\frac{1-j}{m}\choose {k_1}}}{{\frac{1}{2}+\frac{1}{m}\choose 1}}\lambda_{n,3}^{k_1}\right)+o\left(\lambda_{n,0}^{-\frac{1}{2}-\frac{1}{m}}\right)\nonumber\\
&=-\frac{1}{{\frac{1}{2}+\frac{1}{m}\choose 1}}\sum_{j=2}^{\lfloor\frac{m}{2}+1\rfloor}c_{j}(a)\lambda_{n,0}^{-\frac{j}{m}}.\label{asy_eq5}
\end{align}

Suppose that 
$\frac{\lambda_{n,1}}{\lambda_{n,1}}=+\lambda_{n,4}+\lambda_{n,6}+\cdots+\lambda_{n,2s}+\lambda_{n,2s+1},$
where $\lambda_{n,2s+1}=o\left(\lambda_{n,0}^{-\frac{s}{m}}\right)$ and
$
\lambda_{n,2t}=e_{t}(a)\lambda_{n,0}^{-\frac{t}{m}},\,\, 2\leq t\leq s<\frac{m+2}{2}$ for some  $e_{t}(a)\in\C$.
Then from \eqref{asy_eq4}
\begin{align}
&\sum_{k=1}^{\infty}{\frac{1}{2}+\frac{1}{m}\choose k}\left(\lambda_{n,4}+\cdots+\lambda_{n,2s}+\lambda_{n,2s+1}\right)^k\nonumber\\
&+\sum_{j=1}^{\lfloor\frac{m}{2}+1\rfloor}c_{j}(a)\lambda_{n,0}^{-\frac{j}{m}}\sum_{k=1}^{\infty}{\frac{1}{2}+\frac{1-j}{m}\choose k}\left(\lambda_{n,4}+\cdots+\lambda_{n,2s}+\lambda_{n,2s+1}\right)^k+o\left(\lambda_{n,0}^{-\frac{1}{2}-\frac{1}{m}}\right)\nonumber\\
&=-\sum_{j=1}^{\lfloor\frac{m}{2}+1\rfloor}c_{j}(a)\lambda_{n,0}^{-\frac{j}{m}}.\nonumber
\end{align}
Hence,
\begin{align}
&\sum_{k=1}^{\infty}{\frac{1}{2}+\frac{1}{m}\choose k}\left(\lambda_{n,4}+\cdots+\lambda_{n,2s}\right)^k\nonumber\\
&+\sum_{j=1}^{\lfloor\frac{m}{2}+1\rfloor}c_{j}(a)\lambda_{n,0}^{-\frac{j}{m}}\sum_{k=1}^{\infty}{\frac{1}{2}+\frac{1-j}{m}\choose k}\left(\lambda_{n,4}+\cdots+\lambda_{n,2s}\right)^k\nonumber\\
&=-\sum_{j=1}^{\lfloor\frac{m}{2}+1\rfloor}c_{j}(a)\lambda_{n,0}^{-\frac{j}{m}}-{\frac{1}{2}+\frac{1}{m}\choose 1}\lambda_{n,2s+1}+o\left(\lambda_{n,0}^{-\frac{s+2}{m}}\right)+o\left(\lambda_{n,0}^{-\frac{1}{2}-\frac{1}{m}}\right).\label{asy_eq11}
\end{align}
Next,
\begin{align}
&\left(\lambda_{n,4}+\cdots+\lambda_{n,2s}+\lambda_{n,2s+1}\right)^k\nonumber\\
&=\left(e_{2}(a)\lambda_{n,0}^{-\frac{2}{m}}+e_{3}(a)\lambda_{n,0}^{-\frac{3}{m}}+\cdots+e_{s}(a)\lambda_{n,0}^{-\frac{s}{m}}+o\left(\lambda_{n,0}^{-\frac{s}{m}}\right)\right)^k\nonumber\\
&=\sum_{k_1=0}^{k}{k\choose k_1}\left(e_{2}(a)\lambda_{n,0}^{-\frac{2}{m}}+e_{3}(a)\lambda_{n,0}^{-\frac{3}{m}}+\cdots+e_{s}(a)\lambda_{n,0}^{-\frac{s}{m}}\right)^{k-k_1}o\left(\lambda_{n,0}^{-\frac{k_1s}{m}}\right)\nonumber\\
&=\left(e_{2}(a)\lambda_{n,0}^{-\frac{2}{m}}+e_{3}(a)\lambda_{n,0}^{-\frac{3}{m}}+\cdots+e_{s}(a)\lambda_{n,0}^{-\frac{s}{m}}\right)^{k}
+o\left(\lambda_{n,0}^{-\frac{s+2}{m}}\right)\nonumber\\
&=\sum_{\substack{i_p\geq 0,\,j_p\not=j_q\,\,\text{if}\,\, p\not=q\\i_1+\cdots+i_t=k}}\frac{k!}{i_1!\cdots i_t!}e_{j_1}(a)^{i_1}e_{j_2}(a)^{i_2}\cdots e_{j_t}(a)^{i_t}\lambda_{n,0}^{-\frac{i_1j_1+\cdots+i_tj_t}{m}}+o\left(\lambda_{n,0}^{-\frac{s+2}{m}}\right).\nonumber
\end{align}
We use this in \eqref{asy_eq11} to see that the left hand side of \eqref{asy_eq11} is a power series in $\lambda_{n,0}^{-\frac{1}{m}}$. 
Then comparing coefficients of $\lambda_{n,0}^{-\frac{j}{m}}$, $1\leq j\leq s$, we have
\begin{equation}\label{asy_eq12}
-c_{j}(a)=\sum_{\substack{|\alpha|=k\\ \alpha\cdot\beta=j}}{\frac{1}{2}+\frac{1}{m}\choose k}\frac{k!}{\alpha!}e(a)^{\alpha}+\sum_{r=1}^{j-2}c_{r}(a)\sum_{\substack{|\alpha|=k\\ \alpha\cdot\beta=j-r}}{\frac{1}{2}+\frac{1-r}{m}\choose k}\frac{k!}{\alpha!}e(a)^{\alpha}.
\end{equation}
Moreover, if $\frac{s+1}{m}\leq \frac{1}{2}+\frac{1}{m}$ (i.\ e., $s+1\leq \frac{m+2}{2}$) then there exists some constant $e_{s+1}(a)\in\C$ such that
\begin{equation}\label{asy_eq14}
\lambda_{n,2s+1}=e_{s+1}(a)\lambda_{n,0}^{-\frac{s+1}{m}}+o\left(\lambda_{n,0}^{-\frac{s+1}{m}}\right).
\end{equation}
Now we let $\lambda_{n,2s+1}=\lambda_{n,2s+2}+\lambda_{n,2s+3}$ where $\lambda_{n,2s+2}=e_{s+1}(a)\lambda_{n,0}^{-\frac{s+1}{m}}$ and $\lambda_{n,2s+3}=o\left(\lambda_{n,0}^{-\frac{s+1}{m}}\right)$.

If $s+1> \frac{m+2}{2}$ then $\lambda_{n,0}^{-\frac{s+1}{m}}$ could be smaller than the error term $o\left(\lambda_{n,0}^{-\frac{1}{2}-\frac{1}{m}}\right)$ in \eqref{asy_eq11}, and hence we cannot deduce existence of $e_{s+1}(a)$ like we do in \eqref{asy_eq14}. This completes proof of Theorem \ref{eigen_asy}.

\begin{remark}
A first few $e_{j}(a)$ are as follows.
\begin{align}
e_{2}(a)&=-\frac{2m}{m+2}\frac{d_{\ell,2}(a)}{d_{\ell,0}(a)},\quad
e_{3}(a)=-\frac{2m}{m+2}\frac{d_{\ell,3}(a)}{d_{\ell,0}(a)},\nonumber\\
e_{4}(a)&=-\frac{2m}{m+2}\frac{d_{\ell,4}(a)}{d_{\ell,0}(a)}+\frac{3m(m-2)}{(m+2)^2}\left(\frac{d_{\ell,2}(a)}{d_{\ell,0}(a)}\right)^2,\nonumber\\
e_{5}(a)&=-\frac{2m}{m+2}\frac{d_{\ell,5}(a)}{d_{\ell,0}(a)}+\frac{4m(m^2-3m-3)}{(m+2)^3}\frac{d_{\ell,2}(a)}{d_{\ell,0}(a)}\frac{d_{\ell,3}(a)}{d_{\ell,0}(a)},\nonumber\\
e_{6}(a)&=-\frac{2m}{m+2}\frac{d_{\ell,6}(a)}{d_{\ell,0}(a)}+\frac{m(m-6)}{(m+2)^2}\left(\frac{d_{\ell,3}(a)}{d_{\ell,0}(a)}\right)^2\nonumber\\
&+\frac{2m(m-6)}{(m+2)^2}\frac{d_{\ell,2}(a)}{d_{\ell,0}(a)}\frac{d_{\ell,4}(a)}{d_{\ell,0}(a)}+\frac{m(m-2)(9m-2)}{3(m+2)^3}\left(\frac{d_{\ell,2}(a)}{d_{\ell,0}(a)}\right)^3.\nonumber
\end{align}
\end{remark}

%%%%%%%%%%%%%%%%%%%%%%%%%% Appendix A %%%%%%%%%%%%%%%%%%%%%%%%
\appendix
\section{Computing $K_{m,j,k}$} \label{sA}
\renewcommand{\theequation}{A.\arabic{equation}}
\renewcommand{\thetheorem}{A.\arabic{theorem}}
\setcounter{theorem}{0}
\setcounter{equation}{0}
%%%%%%%%%%%%%%%%%%%%%%%%%%%%%%%%%%%%%%%%%%%%%%%%%%%%%%%%%%%%%%
\begin{theorem}
Let $m\geq 3$ be an integer. Then 
\begin{equation}\nonumber
K_m=K_{m,0}=\int_0^{\infty}\left(\sqrt{1+t^m}-t^{\frac{m}{2}}\right)\,dt=\frac{\sqrt{\pi}}{2\cos\left(\frac{\pi}{m}\right)}\frac{\Gamma\left(1+\frac{1}{m}\right)}{\Gamma\left(\frac{3}{2}+\frac{1}{m}\right)}.
\end{equation}
\end{theorem}
\begin{proof}
Substitute $\sqrt{u}=\sqrt{1+t^m}-t^{\frac{m}{2}}$. Then 
\begin{align}
\int_0^{\infty}\left(\sqrt{1+t^m}-t^{\frac{m}{2}}\right)\,dt&=\frac{1}{2^{\frac{2}{m}}m}\int_0^1\left((1-u)^{\frac{2}{m}-1}u^{\frac{1}{2}-\frac{1}{m}-1}+(1-u)^{\frac{2}{m}-1}u^{\frac{3}{2}-\frac{1}{m}-1}\right)\,du\nonumber\\
&=\frac{1}{2^{\frac{2}{m}}m}\left(B\left(\frac{2}{m},\frac{1}{2}-\frac{1}{m}\right)+B\left(\frac{2}{m},\frac{3}{2}-\frac{1}{m}\right)\right),\nonumber
\end{align}
where $B(z,w)$ is the beta function. Then we use the following to complete the proof.
\begin{align}
&\Gamma(z+1)=z\Gamma(z),\quad \Gamma(z)\Gamma(1-z)=-z\Gamma(-z)\Gamma(z)=\frac{\pi}{\sin(\pi z)}\nonumber\\
&B(z,w)=\int_0^1(1-u)^{z-1}u^{w-1}\,du=\frac{\Gamma(z)\Gamma(w)}{\Gamma(z+w)},\quad \Gamma(2z)=\frac{2^{2z-\frac{1}{2}}}{\sqrt{2\pi}}\Gamma(z)\Gamma\left(z+\frac{1}{2}\right).\label{gamma_eq}
\end{align}
\end{proof}

\begin{theorem}
Let $m\geq 3$ and $1\leq k\leq j\leq \frac{m+2}{2}$. Then
\begin{align}\nonumber
K_{m,j,k}&=\int_0^{\infty}\left(\frac{t^{mk-j}}{(1+t^m)^{k-\frac{1}{2}}}-t^{\frac{m}{2}-j}\right)\,dt\nonumber\\
&=\left\{
                    \begin{array}{cl}
-\frac{2}{m}
\quad &\text{if $j=k=1$},\\
&\\
                  -\frac{2k-1}{m+2-2j}B\left(k-\frac{j-1}{m},\,\frac{1}{2}+\frac{j-1}{m}\right)   \quad &\text{if $1\leq k\leq j\leq\frac{m+1}{2}$, $j\not=1$},\\
&\\
                 \frac{2}{m}\left(\ln 2-\frac{1}{1}-\frac{1}{3}-\dots-\frac{1}{2k-5}-\frac{1}{2k-3}\right)  \quad &\text{if $m$ is even, $1\leq k\leq j=\frac{m+2}{2}$.}
                    \end{array}\right.\nonumber
\end{align}
\end{theorem}
\begin{proof}
The case when $j=k=1$ is an easy consequence of 
$$\frac{d}{dt}\left(\sqrt{1+t^m}-t^{\frac{m}{2}}\right)=\frac{m}{2}\left(\frac{t^{m-1}}{(1+t^m)^{\frac{1}{2}}}-t^{\frac{m}{2}-1}\right).$$

Suppose that $1\leq k\leq j\leq\frac{m+1}{2}$, $j\not=1$. Then since
\begin{align}
&\frac{d}{dt}\left(\frac{t^{mk-(j-1)}}{(1+t^m)^{k-\frac{1}{2}}}-t^{\frac{m}{2}-(j-1)}\right)\nonumber\\
&=(mk-(j-1))\left(\frac{t^{mk-j}}{(1+t^m)^{k-\frac{1}{2}}}-t^{\frac{m}{2}-j}\right)-m(k-\frac{1}{2})\left(\frac{t^{m(k+1)-j}}{(1+t^m)^{(k+1)-\frac{1}{2}}}-t^{\frac{m}{2}-j}\right),\nonumber
\end{align}
we have
\begin{align}
\int_0^{\infty}\left(\frac{t^{mk-j}}{(1+t^m)^{k-\frac{1}{2}}}-t^{\frac{m}{2}-j}\right)dt
&=\frac{m(k-1)-(j-1)}{m(k-1)-\frac{m}{2}}\int_0^{\infty}\left(\frac{t^{m(k-1)-j}}{(1+t^m)^{(k-1)-\frac{1}{2}}}-t^{\frac{m}{2}-j}\right)dt\nonumber\\
&=\frac{\Gamma\left(k-\frac{j-1}{m}\right)\Gamma\left(1-\frac{1}{2}\right)}{\Gamma\left(k-\frac{1}{2}\right)\Gamma\left(1-\frac{j-1}{m}\right)}\int_0^{\infty}\left(\frac{t^{m-j}}{(1+t^m)^{\frac{1}{2}}}-t^{\frac{m}{2}-j}\right)dt.\nonumber
\end{align}
Next, we use the substitution $\sqrt{u}=\sqrt{1+t^m}-t^{\frac{m}{2}}$ to show 
$$\int_0^{\infty}\left(\frac{t^{m-j}}{(1+t^m)^{\frac{1}{2}}}-t^{\frac{m}{2}-j}\right)dt=-\frac{2^{\frac{2(j-1)}{m}}}{m}B\left(1-\frac{2(j-1)}{m},\frac{1}{2}+\frac{(j-1)}{m}\right).$$
Finally, we use equations in \eqref{gamma_eq} to complete the proof for $1\leq k\leq j\leq\frac{m+1}{2}$, $j\not=1$.

Finally, if $m$ is even and $j=\frac{m+2}{2}$, then we use integration by parts, for  $R>0$,
\begin{equation}\nonumber
\int_0^{R}\frac{t^{mk-\frac{m}{2}-1}}{\left(t^m+1\right)^{k-\frac{1}{2}}}dt=\left.\frac{1}{m\left(-k+\frac{3}{2}\right)}\frac{t^{m(k-1)-\frac{m}{2}}}{\left(t^m+1\right)^{(k-1)-\frac{1}{2}}}\right|_0^R+\int_0^Rt^{m(k-1)-\frac{m}{2}-1}\frac{1}{\left(t^m+1\right)^{(k-1)-\frac{1}{2}}}dt,
\end{equation}
and hence
$$\int_0^{\infty}\left(\frac{t^{mk-\frac{m}{2}-1}}{\left(t^m+1\right)^{k-\frac{1}{2}}}-\frac{t^{m(k-1)-\frac{m}{2}-1}}{\left(t^m+1\right)^{(k-1)-\frac{1}{2}}}\right)\,dt=-\frac{2}{m\left(2k-3\right)}.$$
Also, one sees that 
$$\int_0^{\infty}\left(\frac{t^{\frac{m}{2}-1}}{\left(t^m+1\right)^{\frac{1}{2}}}-\frac{1}{t+1}\right)\,dt=\frac{2\ln 2}{m}.
$$

\end{proof}
\subsection*{{\bf Acknowledgments}}
The author thanks Mark Ashbaugh, Fritz Gesztesy, Richard Laugesen, Boris Mityagin, Grigori Rozenblioum and Alexander Turbiner for helpful discussions and references. He also thanks Richard Laugesen for reading a part of this manuscript and suggestions for improving its presentation.

{\sc email contact:}  kcshin@math.missouri.edu
\end{document}